\documentclass[11pt]{article}
\usepackage{nsd}

\title{Some bounds on the number of integer points in a ball or on a sphere}
\author{Jiucheng Dai\\
Cornell University \and Rishi Gujjar\\
Cornell University \and Noah Stephens-Davidowitz\thanks{
    \texttt{noahsd@gmail.com}. Supported by NSF Grant No. CCF-2312296, a Packard Foundation Fellowship, and a generous gift from Google.
    }\\
Cornell University}
\date{November 29, 2026}

\begin{document}

\maketitle

\begin{abstract}
    We present multiple upper and lower bounds on the number of integer lattice points in a ball or on a sphere in $d$ dimensions, specifically on the functions
    \[
        r_d(m) := |\{ \vec{z} \in \Z^d \ : \ \|\vec{z}\|^2 = m\}|
    \]
    and
    \[
        N(d,m;\vec{t}) := |\{ \vec{z} \in \Z^d \ : \ \|\vec{z} - \vec{t}\|^2 \leq m\}|
        \; 
    \]
    for $\vec{t} := (t_1,\ldots, t_d) \in \R^d$.

    First, we present relatively simple upper and lower bounds on $r_d(m)$ of the form
    \[
        \binom{m + \alpha d - 1}{m}
        \; 
    \]
    for integers $m$ and some $\alpha > 0$, as well as slightly more complicated generalizations that are \emph{quite} tight. These bounds are proven by a simple technique of showing that if such a bound holds for some small values of $d$, then it must hold for larger values of $d$. We then observe that such bounds hold for small values of $d$ using Jacobi's celebrated exact formulas for $r_d(m)$ for small even values of $d$.
    
    Second, for $m \geq 1$ we prove that
    \[
        N(d,m;\vec{t}) = e^{-\chi(m;\vec{t})} \inf_{\tau > 0} e^{\pi \tau m} \Theta(\tau; \vec{t})
        \; ,
    \]
    where $\Theta(\tau;\vec{t})$ is the $d$-dimensional theta function and the error term $\chi(m;\vec{t}) $ satisfies $ 0 \leq \chi(m;\vec{t}) \leq C\min\{\sqrt{d}, \sqrt{m} \log(d/m+2)\}$ for a small universal constant $C > 0$. 
    We prove the latter result by revisiting Mazo and Odlyzko's beautiful proof that $0 \leq \chi(\alpha d; \vec{t}) \leq C_{\alpha,\delta} \sqrt{d}$ for an unspecified constant $C_{\alpha, \delta} > 0$, where $\delta > 0$ is such that $\|\vec{t}\|^2 \leq (\alpha - \delta) d$. As a consequence, we resolve an open question of Mazo and Odlyzko.
\end{abstract}

\mytableofcontents

\section{Introduction}

We study the functions
\[
    r_d(m) := |\{ \vec{z} \in \Z^d \ : \ \|\vec{z}\|^2 = m\}|
    \; ,
\]
and
\[
    N(d,m) := |\{ \vec{z} \in \Z^d \ : \ \|\vec{z}\|^2 \leq m\}|
    \; ,
\]
where the dimension $d$ and the squared radius $m$ are positive integers. In other words, we are interested in the number of $d$-dimensional vectors with integer coordinates on a sphere or in a ball of radius $\sqrt{m}$.%
\footnote{Our choice to use $d$ for the dimension and $m$ for the square of the radius is not standard, but it avoids a notational issue in the literature where some authors use $n$ for the dimension and $k$ for the squared norm, while others use $k$ for the dimension and $n$ for the squared norm.}
We also study
\[
    N(d,m; \vec{t}) := |\{ \vec{z} \in \Z^d \ : \ \|\vec{z} -\vec{t}\|^2 \leq m\}|
    \; ,
\]
i.e., the number of integer points in such a ball centered around some $\vec{t} \in \R^d$.
In particular, $N(d,m)$ is simply the important special case of $N(d,m;\vec{t})$ when $\vec{t} = \vec0$.\footnote{When considering $N(d,m;\vec{t})$ for $\vec{t} \neq \vec0$, it no longer makes much sense to restrict to the case of integer $m$. But, we ignore this in the introduction.}
All three of these functions of course arise naturally in a wide number of contexts---from number theory to physics to computer science.

In this note, we prove multiple different upper and lower bounds on $N(d,m)$, $N(d,m;\vec{t})$, and $r_d(m)$. These bounds vary in how tight they are in different regimes, their simplicity, and also their novelty. We are particularly interested in simple explicit bounds that are valid for a wide range of parameters $d$ and $m$, particularly the most difficult regime in which $m \approx d$. (In contrast, much literature is interested in obtaining bounds on $N(d,m)$ or $r_d(m)$ that are very tight asymptotically when $d$ is fixed and $m$ goes to infinity, but not in other cases. See, e.g., \cite{ivicLatticePointsLarge2004}.) We hope in particular that this note will serve as a useful reference for researchers who require such bounds on $N(d,m)$, $N(d,m;\vec{t})$, and $r_d(m)$. 

\subsection{The bounds}

The bounds in this note are grouped into three different categories. We present them roughly from the simplest to the most complex. The simplest bounds that we present are arguably trivial (and certainly not novel); we include them only for completeness.

\subsubsection{Basic volume-based bounds} It is of course very natural to guess that the number of integer points in a nice enough subset of $\R^d$ should be roughly equal to the volume of this subset.  So, for completeness, we first include a simple and very well-known proof that $N(d,m;\vec{t})$ is at least the volume of a ball with radius $\sqrt{m} - \sqrt{d}/2$ and at most the volume of a ball with radius $\sqrt{m} + \sqrt{d}/2$. See \cref{thm:volume_estimate}. In fact, this is a special case of a more general (and also well-known) result that applies for arbitrary norms and arbitrary lattices.

Notice that these bounds are quite accurate when $m$ is very large compared to $d$. In particular, for $m > d/4$, such bounds are evidently tight up to a factor of the ratio of the two volumes, which is
\[
    \frac{(1+\sqrt{d/(4m)})^d}{(1-\sqrt{d/(4m)})^d} \approx e^{\sqrt{d^3/m}}
    \; 
\]
for $m \gg d$.
However, for small $m \lesssim d$, volume-based estimates are quite poor. E.g., if $m < d/(2\pi e)$, then the volume of a ball with radius $\sqrt{m}$ is less than one, but $N(d,m;\vec{t})$ can be quite large, even exponential in $d$.

\subsubsection{Bounds based on recurrences and binomial coefficients}

Our second set of bounds are combinatorial and most naturally viewed as bounds on $r_d(m)$,
i.e., the function that counts the number of integer points on a \emph{sphere} rather than in a ball. (This is sometimes called the sum of squares function.) Bounds on $r_d(m)$ can of course be easily translated into bounds on $N(d,m)$ via the simple identity $N(d,m) = \sum_{i=0}^m r_d(i)$. 

Specifically, we prove various upper and lower bounds on $r_d(m)$ of the form
\begin{equation}
    \label{eq:binomial_bound_form_intro}
    c^{d} k^m \binom{m + \alpha d -1}{m}
    \; ,
\end{equation}
for constants $c,k,\alpha > 0$ (and slightly more general bounds).
 Bounds of this form are relatively natural, since we have the simple upper bound
\[
    r_d(m) \leq |\{ \vec{z} \in \Z^d \ : \ \|\vec{z}\|_1 = m\}| \leq 2^{d} \cdot |\{ \vec{z} \in \Z_{\geq 0}^d \ : \ \|\vec{z}\|_1 = m\}| = 2^{d} \cdot \binom{d+m-1}{m}
    \; ,
\]
where $\|\vec{z}\|_1$ is the $\ell_1$ norm (and where the first inequality follows from the injective mapping that maps the vector $\vec{z} = (z_1,\ldots, z_d)$ to $ (\mathrm{sign}(z_1) z_1^2, \ldots, \mathrm{sign}(z_d) z_d^2)$). And, for $m \leq d$, we have the simple lower bound
\[ 
    r_d(m) \geq |\{ \vec{z} \in \{-1,0,1\}^d \ : \ \|\vec{z}\|^2  = m\}| = 2^m \binom{d}{m}
    \; .
\]
So, one might expect that more careful arguments could yield versions of these simple bounds that are reasonably tight for all or most values of $m$ and $d$.

Indeed, we show that bounds of this form follow from  the simple fact that $r_d(m)$ satisfies a nice recurrence relation (\cref{eq:r_d_recurrence}), which is also satisfied by the function $c^d k^m \binom{m+\alpha d-1}{m}$. We use this and a simple induction argument to show that any (upper or lower) bound on $r_{d_i}(m)$ of the form taken in \cref{eq:binomial_bound_form_intro}
for some dimensions $d_1,d_2$ and any positive constants $c$, $k$, and  $\alpha$ yields the same bound for any $d := a d_1 + bd_2$ for positive integers $a$ and $b$.  In fact, our result is a bit more general. See \cref{thm:r_d_inequality_recurrence}.

We then observe that Jacobi's four-, six-, and eight-square theorems (which give precise formulas for $r_4(m)$, $r_6(m)$, and $r_8(m)$ in terms of the divisors of $m$) can be used to obtain rather tight bounds of this form for small $d$. Plugging these bounds into our meta-theorem yields a number of bounds that hold for all $m$ and all but a (short) finite list of dimensions $d$.

For example, our simplest such bounds are
\begin{equation}
    \label{eq:binomial_intro}
    \binom{m+d/2-1}{m} \leq r_d(m) \leq \binom{m+2d-1}{m}
    \; ,
\end{equation}
where the upper bound is valid for all $d$ and the lower bound is valid for $d \geq 5$.%
\footnote{For $d \leq 4$, $r_d(m)$ fluctuates quite a bit with $m$. In particular, $r_1(m)$, $r_2(m)$, and $r_3(m)$ can be zero for arbitrarily large $m$ (e.g., if $m = 8k+7$). $r_4(m)$ is never zero (by Lagrange's celebrated four-square theorem), but we still have that $r_4(2^k) = 24$ for $k \geq 1$ and $r_4(m) \geq 8m$ for any odd $m$. So, the function $r_4(m)$ still fluctuates quite a bit and does not admit any lower bound that is unbounded as $m \to \infty$. It is therefore often necessary to restrict our attention to $d \geq 5$.}
In fact, the upper bound in \cref{eq:binomial_intro} follows trivially from the recurrence relation together with the simple bound $r_1(m) \leq m+1 = \binom{m+1}{1}$, while the lower bound requires Jacobi's theorems.

Notice that these bounds are reasonably tight, since we expect $r_d(m)$ to be very roughly proportional to $m^{\Theta(d)}$  as $m \to \infty$ for fixed  $d > 4$.\footnote{Specifically, for large $m$ in fixed dimension $d > 4$ we expect $r_d(m) \approx \vol(B^d(\sqrt{m})) - \vol(B^d(\sqrt{m-1})) \approx C_d m^{d/2-1} \approx C_d' \binom{m+d/2-1}{m}$, which compares quite nicely with the lower bound in \cref{eq:binomial_intro}. In particular, this must be true on average over $m$, since $\sum_{i=1}^m r_d(i) \approx \vol(B^d(\sqrt{m}))$, as noted above and shown formally in \cref{sec:volume}.} For small $m$ and large $d$, we expect $r_d(m)$ to be roughly equal to $|\{\vec{z} \in \{-1,0,1\}^d \ : \ \|\vec{z}\|^2 = m \}| = 2^m \binom{d}{m} = d^{\Theta(m)}$ as $d \to \infty$ for fixed $m$. The binomial coefficients in \cref{eq:binomial_intro} interpolate nicely between these two regimes, having size $d^{\Theta(m)}$ as $d \to \infty$ for fixed $m$ and $m^{\Theta(d)}$ as $m \to \infty$ for fixed $d$.

Similar ideas were used in \cite{rouseExplicitBoundsSumOfSquares} and \cite[Appendix A]{regevSimpleProofReverse2023}. In particular, the lower bound in \cref{eq:binomial_intro} was originally proven in \cite[Equation (4)]{regevSimpleProofReverse2023}, though only for the (easier) special case when $d \geq 6$ is even. \cite{regevSimpleProofReverse2023} also prove a similar but incomparable upper bound (in a more general setting). In fact, the authors of the present work were inspired by \cite{regevSimpleProofReverse2023} to consider such bounds. 

We illustrate the generality of our approach by proving, e.g., 
\[
r_d(m) \leq \min\{ 1.9^d, 1.04^m \cdot 1.67^d \} \cdot   {m + d/2 -1 \choose m}
\]
for even $d \notin \{2,4,10\}$, and several similar bounds.
(Compare this with the lower bound in \cref{eq:binomial_intro}. This is a slightly weaker and simplified version of \cref{thm:exponential_binomial}, which also includes a lower bound with a similar shape that is quite tight for $m \gtrsim d$.)
Though this bound is loose by a factor of $2^{\Theta(d)}$, the constants are small enough that it is quite tight for a very wide range of parameters. See \cref{sec:some_plots} for plots.

Of course, one can go on to prove innumerable such bounds that will vary in their simplicity, generality, and tightness in different regimes. We hope that our meta-theorem is sufficiently flexible to allow future mathematicians to easily prove whatever bounds are needed in their use cases. However, we resist the urge to enumerate too many bounds of this form, since our next set of bounds is quite tight (and the main benefit of these binomial-based bounds is their relative simplicity).

\subsubsection{Mazo-Odlyzko-style bounds based on the theta function}
Finally, we present bounds on $N(d,m;\vec{t})$ based on the theta function, defined as
\[
    \Theta(\tau; \vec{t}) := \sum_{\vec{z} \in \Z^d} q^{\|\vec{z}- \vec{t}\|^2}
    \; ,
\]
where $q := e^{-\pi \tau}$ for $\tau > 0$.\footnote{It is more common to define $q = e^{\pi i \tau}$ and to take any complex $\tau$ in the upper-half plane, but our parameterization with $\tau$ real is far more convenient for our purposes. The $\pi$ in the exponent is simply a (good) convention, which is convenient, e.g., when applying the Poisson summation formula.} Notice that for $\vec{t} := (t_1,\ldots, t_d) \in \R^d$, we can write
\[
    \Theta(\tau;\vec{t}) = \prod_{i=1}^d \Theta(\tau;t_i)
    \; ,
\]
which in particular means that $\Theta(\tau;\vec{t})$ can be computed efficiently to high accuracy.\footnote{In more detail, if, say, $\tau \geq 1$, then the one-dimensional series $\Theta(\tau;t)$ converges very rapidly, with only $O(\sqrt{n})$ terms needed to compute the series to within a factor of $1+2^{-n}$. For $\tau < 1$, one can use the Poisson summation formula to write $\Theta(\tau;t)$ as a similarly rapidly converging series.}
Furthermore, Markov's inequality immediately implies that
\begin{equation}
    \label{eq:theta_upper_bound_intro}
    N(d,m;\vec{t}) \leq \inf_{\tau > 0} q^{-m} \Theta(\tau; \vec{t})
    \; .
\end{equation}
(Though the right-hand side of this equation is an infimum of a $d$-dimensional series, it can still be computed to high accuracy quite efficiently by again noting that $\Theta(\tau;\vec{t})$ can be written as a product of $d$ rapidly converging one-dimensional series. We discuss this more in \cref{sec:finite_sum_defs} and in the appendix.)

In a beautiful paper, Mazo and Odlyzko~\cite{mazoLatticePointsHighdimensional1990} proved that \cref{eq:theta_upper_bound_intro} is quite tight if we fix $m := \alpha d$ for some constant $\alpha > 0$ and take $d \to \infty$. Specifically, they proved that for any constants $0 < \delta < \alpha$ and $\vec{t} \in \R^d$ with $\|\vec{t}\|^2 \leq (\alpha-\delta) d$, there exists $C_{\alpha, \delta} > 0$ such that
\[
    e^{-C_{\alpha, \delta} \sqrt{d}} \inf_{\tau > 0}  \frac{\Theta(\tau; \vec{t})}{q^{\alpha d}} \leq N(d,\alpha d;\vec{t}) \leq \inf_{\tau > 0} \frac{\Theta(\tau; \vec{t})}{q^{\alpha d}}
    \; .
\]
Their proof is quite elegant, but it (1) only works for $m = \alpha d$ for constant $\alpha$ (and for $\|\vec{t}\|^2$ bounded away from $\alpha d$), (2) does not give an explicit bound on the constant $C_{\alpha, \delta}$, and (3) does not provide much intuition. (At a technical level, all of these ``issues'' arise from the elegant compactness arguments that Mazo and Odlyzko use in their proof.)

We observe that Mazo and Odlyzko's celebrated result is actually perhaps best viewed as simply a concentration bound on the \emph{discrete Gaussian distribution}, $D(\tau;\vec{t})$, i.e., the distribution over $\Z^d - \vec{t}$ defined by
\[
    \Pr_{\vec{X} \sim D(\tau;\vec{t})}[\vec{X} = \vec{x}] := q^{\|\vec{x}\|^2}/\Theta(\tau;\vec{t})
\]
for any $\vec{x} \in \Z^d - \vec{t}$. Notice in particular that \cref{eq:theta_upper_bound_intro} will be relatively tight if there exists a $\tau > 0$ such that $D(\tau;\vec{t})$ is relatively concentrated on vectors $\vec{x} \in \Z^d - \vec{t}$ with $m-\eps \leq \|\vec{x}\|^2 \leq m$ for some small $\eps$. 

We therefore prove an explicit and more general bound than \cite{mazoLatticePointsHighdimensional1990} by directly proving concentration of this distribution. Our methods are quite elementary, though the proof is a bit delicate in places. Specifically, we carefully apply bounds on $\Theta$ series and related series to just directly bound the variance of $\|\vec{X}\|^2$ when $\vec{X} \sim D(\tau;\vec{t})$. Eventually, we prove that 
\[
    e^{-3 - \min\{4\sqrt{d}, 8 \sqrt{m}\log'(d/m)\}} \inf_{\tau > 0}  \frac{\Theta(\tau; \vec{t})}{q^{m}} \leq N(d,m;\vec{t}) \leq \inf_{\tau > 0} \frac{\Theta(\tau; \vec{t})}{q^m}
    \; ,
\]
where we write $\log'(x) := \max\{1,\log x\}$. (See \cref{thm:main_theta_theorem}.)
In particular, just the fact that our bound holds for any $\vec{t}$ when $m = d/4$ resolves an open question asked in \cite[Section 3]{mazoLatticePointsHighdimensional1990}.\footnote{In our notation, they asked whether $\vec{t} := (1/2,\ldots, 1/2) \in \R^d$ satisfies $N(d,d/4;\vec{t}) \leq e^{O(\sqrt{d})} N(d,d/4;\vec{t}')$ for all $\vec{t}' \in \R^d$. They were able to prove in \cite[Theorem 2]{mazoLatticePointsHighdimensional1990} the analogous result for $m = \alpha d$ for any constant $\alpha > 1/4$ and speculated that this also holds for $\alpha = 1/4$. Our result shows that this does in fact hold for $\alpha = 1/4$, once one observes that for any $\tau$, $\Theta(\tau;\vec{t})$ is minimized at $\vec{t} = (1/2,\ldots, 1/2)$ (as was observed in \cite{mazoLatticePointsHighdimensional1990} using the Jacobi triple product formula).} (To be clear, we made little attempt to optimize the constants in the lower bound, and it is clear from the proof that they are not tight.)

Notice that this bound applies for any $\vec{t} \in \R^d$ and any $m \geq 1$,
is effective, and quite tight.  
For example, when $m$ is a constant, \cref{thm:main_theta_theorem} gives bounds on $N(d,m;\vec{t})$ that are tight up to a factor that is polynomial in $d$. This is clearly the best possible for any continuous function, since, e.g., $N(d,m;\vec{0})$ jumps by a $\poly(d)$ factor at each positive constant integer $m$. In general, we do not know whether an error term of the form $e^{-O(\min\{\sqrt{d}, \sqrt{m} \log'(d/m) \})}$ is the best possible. But, it seems to be essentially the best possible using our methods. So, any significant improvement will need to use substantially different ideas.

(See also recent work of Niksi{\'n}ski and Wr{\'o}bel~\cite{niksinskiDimensionfreeEstimatesDiscrete2026}, which showed quantitative estimates on $r_d(m)$ and $N(d,m;\vec0)$ for $m \leq O(d)$ using the theta function. And, see \cite{elkiesPackingDensitiesSuperballs1991,aggarwalGapSETHHardness2018} for analogues of \cite{mazoLatticePointsHighdimensional1990} for other norms.)

\subsection{AI disclosure}

AI tools were used for (1) optimizing the Mathematica code used to produce the plots in \cref{sec:some_plots}; (2) a literature search to try to find earlier resolution of the open problem in \cite{mazoLatticePointsHighdimensional1990} (without success); and (3) proofreading. The tools found typos and minor errors in earlier versions of this work.

\subsection*{Acknowledgements}

    The authors thank Huck Bennett, Alexander Golovnev, Robert Kleinberg, and Oded Regev for helpful discussions.

\section{Preliminaries}

\subsection{Binomial coefficients}

For $x, k \in \R$, we define the generalized binomial coefficient as 
\[
    \binom{x}{k} := \frac{\Gamma(x+1)}{\Gamma(k+1)\Gamma(x-k+1)}
    \; ,
\]
where $\Gamma$ is the Gamma function, and we naturally interpret $1/\Gamma(-n) = 0$ for non-negative integers $n$, $\Gamma(x+1)/\Gamma(x-n+1) = \prod_{i=1}^n (x-i+1)$ for non-negative integers $n$ (even if $x$ is a negative integer, for which $\Gamma(x+1)$ is undefined), and $\binom{x}{k} = 0$ for negative integers $k$. In particular, for non-negative integers $x$ and $k$, this is the familiar binomial coefficient. We will mostly be interested in cases when $x > 0$ and $x, k \in \tfrac{1}{2} \Z$.

We will use the following variant of a well-known combinatorial identity for the binomial coefficients.

\begin{lemma}
    \label{lem:2dim_hockey_stick}
    For any non-negative integers $m$, $k_1$, and $k_2$ and any $x_1,x_2 \in \R$, we have 
    \[
        \sum_{i=0}^{m}
            {i+x_1-1 \choose i-k_1}
            {m-i+x_2-1 \choose m-i-k_2}=
            {m+x_1+x_2-1 \choose m-k_1-k_2}
        \; .
    \]
\end{lemma}

\begin{proof}
    We work with formal power series. We first recall that for any $x \in \R$,
    \begin{equation}
        \label{eq:generalized_binomial_gf}
        (1-z)^{-x}
        =
        \sum_{j=0}^{\infty} {j+x-1 \choose j} z^j
        \; .
    \end{equation}
    Indeed, this is just the generalized binomial theorem.
    Equivalently, \cref{eq:generalized_binomial_gf} may be viewed purely as an identity of formal power series.

    For any non-negative integer $k$, \cref{eq:generalized_binomial_gf} gives
    \[
        \sum_{i=0}^{\infty} {i+x-1 \choose i-k} z^i
        =
        z^k \sum_{j=0}^{\infty} {j+k+x-1 \choose j} z^j
        =
        z^k(1-z)^{-(x+k)}
        \; ,
    \]
    where the first equality follows from the change of variables $j=i-k$.

    Applying this identity with $(x,k)=(x_1,k_1)$ and $(x,k)=(x_2,k_2)$, multiplying, and comparing the coefficient of $z^m$, we get
    \begin{align*}
        \sum_{i=0}^{m}
            {i+x_1-1 \choose i-k_1}
            {m-i+x_2-1 \choose m-i-k_2}
        &=
        [z^m]\,
        z^{k_1+k_2}(1-z)^{-(x_1+x_2+k_1+k_2)}
        \\
        &=
        [z^{m-k_1-k_2}]\,
        (1-z)^{-(x_1+x_2+k_1+k_2)}
        \\
        &=
        {m+x_1+x_2-1 \choose m-k_1-k_2}
        \; ,
    \end{align*}
    again by \cref{eq:generalized_binomial_gf}.
\end{proof}

\subsection{Simple bounds on \texorpdfstring{$r_d$}{the sum of squares function} for small \texorpdfstring{$d$}{d}}
\label{sec:simple_bounds_prelims}

Recall that $r_d(m) := |\{\vec{z} \in \Z^d \ : \ \|\vec{z}\|^2 = m\}|$. We first derive some simple bounds on $r_6(m)$ and $r_8(m)$ which we will use extensively to derive results for small dimensions.

\begin{lemma}
    \label{r_6 bound}
    For any integer $m \geq 1$, we have
    \[ 10.2 m^2 < C_1 m^2 < r_6(m) < C_2 m^2 < 21.13 m^2 
    \; ,\]
    where 
    \[
        C_1 := 12 \prod_{\substack{p \equiv 3 \bmod 4\\ p \text{ prime}}} (1-1/p^2) \approx 10.3
    \]
    and
    \[
        C_2 := 20 \prod_{\substack{p \equiv 1 \bmod 4\\ p \text{ prime}}} \frac{1}{1-1/p^2} \approx 21.1
        \; .
    \]
    
    In fact,
        \[
         \liminf_{m \to \infty} \frac{r_6(m)}{m^2} = C_1
         \;,
    \]
    and
    \[
        \limsup_{m \to \infty} \frac{r_6(m)}{m^2} = C_2
       \; .
    \]
\end{lemma}
\begin{proof}
    We recall Jacobi's six-square theorem which states that
    \begin{align*} r_6(m) 
        &= 4m^2 \sum_{d \mid m} \frac{1}{d^2} \cdot (4\chi(d) -\chi(m/d)) \\
        &= 4 m^2(4-\chi(m)) + 4 m^2 \sum_{\substack{d \mid m\\d > 1}} \frac{1}{d^2} \cdot (4\chi(d) -\chi(m/d))
        \; .
    \end{align*}
    where $\chi(k) := \left( \frac{-4}{k} \right) \in \{-1,0,1\}$ is the Kronecker symbol satisfying
    \[
        \chi(k) = \begin{cases}
            1 & k \equiv 1 \bmod 4\\
            -1 & k \equiv 3 \bmod 4\\
            0 & k \equiv 0 \bmod 2 \; .
        \end{cases}
    \]
    (See, e.g.,~\cite{Grosswald1985}.)

    The result then follows by partitioning into cases according to $m \bmod 4$.

    The easiest case is when $m$ is even. Then,
    we see that
    \[
        r_6(m) = 16m^2 + 16m^2 \sum_{\substack{d \mid m \\ d \equiv 1 \bmod 4\\d > 1}} \frac{1}{d^2}- 16m^2 \sum_{\substack{d \mid m \\ d \equiv 3 \bmod 4}} \frac{1}{d^2} - 4m^2 \sum_{\substack{d \mid m \\ m/d \equiv 1 \bmod 4}}  \frac{1}{d^2} + 4m^2 \sum_{\substack{d \mid m \\ m/d \equiv 3 \bmod 4}}  \frac{1}{d^2} 
        \; .
    \]
    The odd divisor sums are clearly at least $-\sum_{k=0}^\infty 1/(4k+3)^2$ and at most $\sum_{k=1}^\infty 1/(4k+1)^2$. The even divisor sums are trivially bounded in absolute value by $\sum_{k=1}^\infty 1/(2k)^2 = \pi^2/24$. A quick numerical check shows that both the upper and lower bounds hold in this case.
    
    If $m \equiv 1 \bmod 4$, then we see that
    \begin{align*}
        r_6(m)
        &= 12m^2 \sum_{\substack{d \mid m \\ d \equiv 1 \bmod 4}} \frac{1}{d^2}- 12m^2 \sum_{\substack{d \mid m \\ d \equiv 3 \bmod 4}} \frac{1}{d^2}
        \; .
    \end{align*}
    The upper bound is then trivial in this case. The lower bound in this case follows by noting that 
    \[
        \sum_{\substack{d \mid m \\ d \equiv 1 \bmod 4}} \frac{1}{d^2}- \sum_{\substack{d \mid m \\ d \equiv 3 \bmod 4}} \frac{1}{d^2} \geq \prod_{\substack{p \mid m\\ p \equiv 3 \bmod 4\\ p \text{ prime}}} (1-1/p^2) > \prod_{\substack{p \equiv 3 \bmod 4\\ p \text{ prime}}} (1-1/p^2)
        \; ,
    \]
    as needed.
    
    If $m \equiv 3 \bmod 4$, then we see that
     \begin{align*}
        r_6(m)
        &= 20m^2 \sum_{\substack{d \mid m \\ d \equiv 1 \bmod 4}} \frac{1}{d^2} - 20m^2 \sum_{\substack{d \mid m \\ d \equiv 3 \bmod 4}} \frac{1}{d^2} 
        \; .
    \end{align*}
    The lower bound is trivial in this case. The upper bound in this case follows by noting that
    \[
      \sum_{\substack{d \mid m \\ d \equiv 1 \bmod 4}} \frac{1}{d^2} - \sum_{\substack{d \mid m \\ d \equiv 3 \bmod 4}} \frac{1}{d^2} \leq \prod_{\substack{p^k \mid m\\ p \equiv 1 \bmod 4\\ \text{$p$ prime, $k$ maximal}}} \sum_{i=0}^k \frac{1}{p^{2i}} < \prod_{\substack{p \equiv 1 \bmod 4\\ p \text{ prime}}} \frac{1}{1-1/p^2}
        \; ,
    \]
    as claimed.

    So, in all cases, the upper and lower bounds hold. We now briefly observe that the limit infimum and limit supremum are as claimed. To see that the limit infimum is as claimed, take $m$ to be the product of the first $2n$ primes that are congruent to $3$ mod $4$, and notice that as $n \to \infty$ the above proof of the lower bound is tight. For the limit supremum, take $m$ to be the product of $p^n$ for the first $n$ primes $p$ congruent to $1$ mod $4$ times, say, the $n$th prime that is congruent to $3$ mod $4$.
\end{proof}

$r_8(m)$ is considerably easier to work with than $r_6(m)$. So, the proof of our bounds for $r_8(m)$ below is considerably easier, and the resulting bounds are more obviously tight.\footnote{In more detail, $r_8(2) = 112 = 14 \cdot 2^3$, so that the lower bound in \cref{r_8 bound} is in fact precisely tight for one value of $m$. It is also asymptotically tight for infinitely many $m$ in the sense that if $m = 2p$ for an odd prime $p$, then $r_8(m) = 14m^3 + 112 \approx 14m^3$. The upper bound is asymptotically tight in the sense that, $r_8(k!)/(k!)^3 = (1+o(1)) 16\zeta(3)$ as $k \to \infty$.}

\begin{lemma}
    \label{r_8 bound}
    For any integer $m \geq 1$, we have 
    \[ 14m^3 \leq r_8(m) < 16 \zeta(3) m^3 < 19.3m^3
    \; .\]
\end{lemma}
\begin{proof}
    We recall Jacobi's eight-square theorem, which states that
    \begin{equation}
        \label{eq:jacobi_eight}
     r_8(m) = 16\sum_{d \mid m}{(-1)^{d+m}d^3}.
    \end{equation}
    (See, e.g.,~\cite{Grosswald1985}.)
    In particular, by making the change of variables
    $k = m/d$, we immediately see that
    \[
        r_8(m) = 16 m^3\sum_{k \mid m} (-1)^{m/k+m} \frac{1}{k^3} <  16m^3 \sum_{k \geq 1} \frac{1}{k^3} = 16 \zeta(3) m^3
        \; ,
    \]
    as needed. 
    
    For the lower bound, we note that if $m$ is odd, all terms in \cref{eq:jacobi_eight} are positive, so that we trivially have $r_8(m) \geq 16m^3$ in this case.
    
    For even $m$, we have
    \begin{align*}
        r_8(m) &= 16\sum_{d \mid m}{(-1)^{d}d^3} \\
        &= 16\sum_{\substack{d \mid m \\ d \text{ even}}}{d^3} - 16\sum_{\substack{d \mid m \\ d \text{ odd}}}{d^3}
        \; .
    \end{align*}
    Noting that every odd divisor is half of an even divisor, this yields
    \begin{align*}
        r_8(m) &= 16\sum_{\substack{d \mid m \\ d \text{ even}}}{d^3} - 16\sum_{\substack{d \mid m \\ d \text{ odd}}}{d^3} \\
        &\geq 16\sum_{\substack{d \mid m \\ d \text{ even}}}{d^3} - 16\sum_{\substack{d \mid m \\ d \text{ even}}}{(d/2)^3} \\
        &\geq 14 m^3
        \; ,
    \end{align*}
    as needed.
\end{proof}
\begin{corollary}
    \label{binomial_r_6_r_8_bounds}
    For non-negative integers $m$, we have 
    \[ {m+2 \choose 2}\leq r_6(m)
    \; ,
    \] 
    and 
    \[
    {m+3 \choose 3} \leq r_8(m) \; .\]
\end{corollary}
\begin{proof}
    For $m \geq 1$, this follows as a consequence of \cref{r_6 bound} and \cref{r_8 bound}. When $m=0$, it is apparent that equality holds.
\end{proof}

We now consider lower bounds on $r_d(m)$ for small odd integers $d$, which are a bit more difficult to prove than the case of even $d$.

\begin{lemma}
    \label{r_5_bound}
    For any integer $m \geq 0$, we have 
    \[ r_5(m) \geq {m+3/2 \choose 3/2}\]
\end{lemma}
\begin{proof}
    We will prove the result for $m \geq 100$; one can directly verify that the bound holds for $m<100$ (e.g., using \cref{eq:r_d_recurrence}).
    
    Recall Jacobi's four-square theorem, which in particular states that for odd $m$ we have
    \[ 
        r_4(m) = 8\sum_{d \mid m}{d} \geq 8m 
        \; .
    \]
    (See, e.g.,~\cite{Grosswald1985}.)
    Let $k := \floor{\sqrt{m}}$.  By applying the simple identity
    \[
    r_{5}(m)=\sum_{i=0}^{m}{r_{1}(i) \cdot r_{4}(m-i)}
    \; 
\]
(see~\cref{eq:r_d_recurrence}), we see that
    \begin{align*}
        r_5(m) &= \sum_{i=0}^{m}{r_1(i) \cdot r_4(m-i)} 
        \\
        &=r_4(m) + 2 \sum_{j=1}^{k}{r_4(m-j^2)} \\
        &\geq 16 \sum_{\substack{1 \leq j \leq k \\ m-j \text{ odd}}} (m-j^2)\\
        &\geq 4 \sum_{j=1}^{k-2} (m-j^2)\\
        &= 4(k-2) (m - k^2/3 + 5k/6 - 1/2 )
        \; ,
    \end{align*}
    where the last inequality uses the fact that for every $1 \leq j \leq k-1$, we have $4(m-j^2) \geq (m-j^2) + (m-(j-1)^2)$. By plugging in $k \geq \sqrt{m}-1$, we see that
    \[
        r_5(m) \geq \frac{2}{3} \cdot (\sqrt{m} - 3)(4m+9\sqrt{m}-10) \geq 2m^{3/2}
        \; ,
    \]
    where in the last step we have used the fact that $m \geq 100$. The result then follows from the simple inequality $\binom{m+3/2}{3/2} \leq (m+3/2)^{3/2} \leq 2m^{3/2}$, where the last step again uses the fact that $m \geq 100$.
\end{proof}

\begin{lemma}
    \label{lem:r_7_and_r_9_lower_bound}
    For any integer $m \geq 0$, we have
    \[r_7(m) \geq {m + 5/2 \choose 5/2}
    \; ,
    \]
    and 
    \[
        r_9(m) \geq {m+7/2 \choose 7/2} \; .
    \]
\end{lemma}
The techniques used in the proof of \cref{lem:r_7_and_r_9_lower_bound} are essentially the same as in \cref{r_5_bound} so the details of the proof have been left to \cref{app:r_7_and_r_9_lower_bounds}.

\subsection{Some bounds on Gaussian sums and integrals}

We will need the following simple bounds on Gaussian sums.

\begin{claim}
    \label{claim:Gaussian_tail_bound}
    For any $\tau > 0$ and $y > 0$, we have
    \[
        \int_{y}^\infty q^{x^2} \dl{x} \leq \frac{q^{y^2}}{2\pi \tau y}
        \; ,
    \]
    \[
        \int_y^\infty x^2 q^{x^2} \dl{x} \leq \Big( \frac{y}{2\pi \tau} + \frac{1}{4 \pi^2 \tau^2 y}\Big)\cdot q^{y^2} 
        \; ,
    \]
    and
    \[
        \int_y^\infty x^4 q^{x^2} \dl{x} \leq \Big( \frac{y^3}{2\pi \tau} + \frac{3 y}{4\pi^2 \tau^2} + \frac{3}{8\pi^3 \tau^3 y}\Big) \cdot q^{y^2}
        \; ,
    \]
    where $q := e^{-\pi \tau}$.
\end{claim}
\begin{proof}
    We have
    \[
        \int_{y}^\infty q^{x^2} \dl{x} \leq \int_{y}^\infty (x/y) \cdot q^{x^2} \dl{x} = \frac{q^{y^2}}{2\pi \tau y}
        \; ,
    \]
    as claimed.

    Furthermore, by integration by parts, we have
    \[
     \int_y^{\infty} x^2 q^{x^2} \dl{x} = \frac{y}{2\pi \tau} q^{y^2} + \frac{1}{2\pi \tau} \int_{y}^\infty q^{x^2} \dl{x} \leq \frac{y}{2\pi \tau} q^{y^2} + \frac{1}{4 \pi^2 \tau^2 y} q^{y^2}
     \; ,
    \]
    as claimed.

    Similarly,
    \[
        \int_{y}^\infty x^4 q^{x^2} \dl{x} = \frac{y^3}{2\pi \tau} q^{y^2} + \frac{3 y}{4\pi^2 \tau^2} q^{y^2} + \frac{3}{4\pi^2 \tau^2} \int_{y}^\infty q^{x^2} \dl{x} \leq  \frac{y^3}{2\pi \tau} q^{y^2} + \frac{3 y}{4\pi^2 \tau^2} q^{y^2} + \frac{3}{8\pi^3 \tau^3 y} q^{y^2}
        \; ,
    \]
    as claimed.
\end{proof}

From the claim above, we obtain bounds on Gaussian sums as well.

\begin{claim}
    \label{clm:Gaussian_sum_bounds}
    For any $\tau > 0$, $t \geq 0$, and integer $k \geq t + \sqrt{2/(\pi \tau)}$, we have
    \[
        \sum_{z=k+1}^\infty q^{(z-t)^2} \leq \frac{q^{w^2}}{2\pi \tau w}
        \; ,
    \]
    \[
        \sum_{z = k+1}^\infty (z-t)^2 q^{(z-t)^2} \leq  \Big( \frac{w}{2\pi \tau} + \frac{1}{4\pi^2 \tau^2 w} \Big) \cdot q^{w^2}
        \; ,
    \]
    and
    \[
        \sum_{z = k+1}^\infty (z-t)^4 q^{(z-t)^2} \leq  \Big( \frac{w^3}{2\pi \tau} + \frac{3w}{4\pi^2 \tau^2} + \frac{3}{8\pi^3 \tau^3 w} \Big) \cdot q^{w^2}
        \;, 
    \]
    where $q := e^{-\pi \tau}$ and $w := k-t$.
\end{claim}
\begin{proof}
    Notice that the function $f_\tau(x) := q^{x^2}$ is decreasing as a function of $x$ for $x \geq 0$. Therefore, for all $z \geq 1+t$, we have
    \[
        q^{(z-t)^2} \leq \int_{z-t-1}^{z-t} q^{x^2} \dl{x}
        \; .
    \]
    Applying this to the summation, we see that
    \[
        \sum_{z=k+1}^\infty q^{(z-t)^2} \leq \int_{w}^\infty q^{x^2}\dl{x} \leq \frac{q^{w^2}}{2\pi \tau w}
        \; ,
    \]
    where the last inequality is \cref{claim:Gaussian_tail_bound}.

    Similarly, the function $g_\tau(x) := x^2 q^{x^2}$ is decreasing as a function of $x$ for $x \geq 1/\sqrt{\pi \tau}$. So, by the same argument, we see that
    \[
        \sum_{z=k+1}^\infty (z-t)^2 q^{(z-t)^2} \leq \int_{w}^\infty x^2 q^{x^2} \dl{x} \leq \Big( \frac{w}{2\pi \tau} + \frac{1}{4\pi^2 \tau^2 w} \Big) \cdot q^{w^2}
        \; ,
    \]
    where the last inequality is again \cref{claim:Gaussian_tail_bound}.

    Finally, noting that $h_\tau(x) := x^4 q^{x^2}$ is decreasing as a function of $x$ for $x \geq \sqrt{2/(\pi \tau)}$, we see that
    \[
        \sum_{z=k+1}^\infty (z-t)^4 q^{(z-t)^2} \leq \int_{w}^\infty x^4 q^{x^2} \dl{x} \leq \Big( \frac{w^3}{2\pi \tau} + \frac{3w}{4\pi^2 \tau^2} + \frac{3}{8\pi^3 \tau^3 w} \Big) \cdot q^{w^2}
        \; ,
    \]
    where again the last inequality is \cref{claim:Gaussian_tail_bound}.
\end{proof}

We will also need the following result from \cite{regevInequalityGaussiansLattices2017}.\footnote{%
Note that the notation in \cite{regevInequalityGaussiansLattices2017} is quite different from ours, and the result that they prove is stronger and more general than what we need. So, we clarify here how their result implies \cref{thm:second_moment_lower_bound}. For a lattice $\lat \subset \R^d$ and vector $\vec{x} \in \R^d$, \cite{regevInequalityGaussiansLattices2017} uses the notation $D_{\lat + \vec{x}}$ for the distribution on the lattice coset $\lat + \vec{x}$ defined by
\[
    \Pr_{\vec{X} \sim D_{\lat + \vec{x}}}[\vec{X} = \vec{y}] := \frac{e^{-\pi \|\vec{y}\|^2}}{\sum_{\vec{w} \in \lat} e^{-\pi \|\vec{w} + \vec{x}\|^2}}
    \; .
\]
By taking the lattice $\lat := \sqrt{\tau} \Z^d$ and $\vec{x} := -\sqrt{\tau} \vec{t}$, one sees that $D_{\lat+\vec{x}}/\sqrt{\tau} = D(\tau;\vec{t})$. Plugging this in  to \cite[Corollary 3.2]{regevInequalityGaussiansLattices2017}, one obtains the first inequality in \cref{thm:second_moment_lower_bound}. The second inequality follows by the trivial observation that the second moment of $\vec{Y} \sim D(\tau;\vec0)$ is at least the second moment of the distribution obtained by conditioning on $\vec{Y} \in \{-1,0,1\}^d$, which is precisely $2dq/(1+2q)$.} Here, for $\tau > 0$ and $\vec{t} \in \R^d$, we write $D(\tau;\vec{t})$ for the probability distribution over $\Z^d - \vec{t}$ defined so that for every
$\vec{x} \in \Z^d - \vec{t}$, we have
\[
    \Pr_{\vec{X} \sim D(\tau; \vec{t})}[\vec{X} = \vec{x}] := \frac{q^{\|\vec{x}\|^2}}{\sum_{\vec{z} \in \Z^d} q^{\|\vec{z} - \vec{t}\|^2}}
    \; ,
\]
where $q := e^{-\pi \tau}$.

\begin{theorem}[{\cite[Corollary 3.2]{regevInequalityGaussiansLattices2017}}]
    \label{thm:second_moment_lower_bound}
    For any $\tau > 0$ and $\vec{t} \in \R^d$, 
    \[
        \expect_{\vec{X} \sim D(\tau;\vec{t})}[\|\vec{X}\|^2] \geq \Big\|\expect_{\vec{X} \sim D(\tau;\vec{t})}[\vec{X}] \Big\|^2 + \expect_{\vec{Y} \sim D(\tau;\vec0)}[\|\vec{Y}\|^2] \geq \frac{2dq}{1+2q}
        \; ,
    \]
    where $q := e^{-\pi \tau}$.
\end{theorem}

\subsection{The Poisson summation formula applied to some Gaussian sums}

Let 
\[\Theta(\tau;t) := \sum_{z \in \Z} e^{-\pi \tau (z-t)^2}
\; ,
\]
\[
    \mu(\tau;t) := \sum_{z \in \Z} (z-t)^2 e^{-\pi \tau (z-t)^2}
    \; ,
\]
and
\[
    M(\tau;t) := \sum_{z \in \Z} (z-t)^4 e^{-\pi \tau (z-t)^2}
    \; .
\]

\begin{fact}
    \label{fact:PSF_theta_mu_M}
    For any $\tau > 0$ and $t \in \R$,
    \[
        \Theta(\tau;t) = \frac{1}{\sqrt{\tau}} + \frac{2}{\sqrt{\tau}} \cdot \sum_{z=1}^\infty \cos(2\pi t z)p^{z^2}
        \; ,
    \]
    \[
        \mu(\tau;t) = \frac{1}{2\pi \tau} \Theta(\tau;t) - \frac{2}{\tau^{5/2}} \sum_{z=1}^\infty \cos(2\pi z t) z^2 p^{z^2}
        \; ,
    \]
    and
    \[
        M(\tau;t) = \frac{3}{\pi \tau} \mu(\tau;t) - \frac{3}{4\pi^2 \tau^2} \Theta(\tau;t) + \frac{2}{\tau^{9/2}} \sum_{z=1}^\infty \cos(2\pi z t) z^4 p^{z^2}
        \; ,
    \]
    where $p := e^{-\pi/\tau}$.
\end{fact}
\begin{proof}
    All of these identities follow immediately from the Poisson summation formula, which tells us that 
    \[
        \sum_{z \in \Z} f(z) = \sum_{z \in \Z} \widehat{f}(z)
        \; ,
    \]
    for sufficiently nice functions $f : \R \to \C$,
    where 
    \[
        \widehat{f}(w) := \int_{-\infty}^\infty e^{-2\pi i w x} f(x) \dl{x}
    \]
    is the Fourier transform of $f$. Specifically, one obtains the identities by applying this formula to the functions $f_0(x) := e^{-\pi \tau (x-t)^2}$, $f_1(x) := (x-t)^2 e^{-\pi \tau (x-t)^2}$, and $f_2(x) := (x-t)^4 e^{-\pi \tau (x-t)^2}$ respectively.
\end{proof}

\section{Simple bounds in terms of the volume of a ball}
\label{sec:volume}

A very natural estimate for $N(d,m;\vec{t})$ is simply 
\begin{equation}
    \label{eq:volume_approx}
    N(d,m;\vec{t}) \approx \vol(B^d(\sqrt{m})) = \pi^{d/2} m^{d/2}/\Gamma(d/2+1) \approx (\pi d)^{-1/2} \cdot (2\pi e m/d)^{d/2}
\; ,
\end{equation}
where $B^d(r)$ is the $d$-dimensional Euclidean ball of radius $r$. 

This approximation is in fact quite accurate for large $m$, as the following theorem shows. Note that this is a simple and well known result (already known to Gauss, at least in the special case when $d = 2$), which we include here for completeness.

\begin{theorem}
    \label{thm:volume_estimate}
    For any $\vec{t} \in \R^d$ and $r \geq 0$, we have
    \[
        \vol(B^d(r-\sqrt{d}/2)) \leq N(d,r^2;\vec{t}) \leq \vol(B^d(r+\sqrt{d}/2))
        \; .
    \]
    (Here we interpret $\vol(B^d(r')) = 0$ for $r' \leq 0$.)
\end{theorem}
\begin{proof}
    Notice that for any $\vec{s}_1, \vec{s}_2$ and any $r^* \geq 0$, we have
    \begin{equation}
        \label{eq:distance_on_torus}
        N(d,(r^*)^2;\vec{s}_1) \leq N(d,(r^* + \delta)^2;\vec{s}_2) \leq N(d,(r^* + \sqrt{d}/2)^2;\vec{s}_2)
        \; ,
    \end{equation}
    where $\delta := \min_{\vec{z} \in \Z^d} \|\vec{s}_1 - \vec{s}_2 - \vec{z}\| = \dist(\vec{s}_1 - \vec{s}_2, \Z^d) \leq \sqrt{d}/2$. 
    So, to prove the theorem, it suffices to prove that for any $r' \geq 0$, 
    \begin{equation}    
        \label{eq:min_vs_average}
        \min_{\vec{s}} N(d,(r')^2;\vec{s}) \leq \vol(B^d(r')) \leq \max_{\vec{s}} N(d,(r')^2;\vec{s})
        \; .
    \end{equation}
    Then, for the upper bound, we can simply set $r^* := r$ and $\vec{s}_1 := \vec{t}$ and then take $\vec{s}_2$ such that $N(d,(r^* + \sqrt{d}/2)^2;\vec{s}_2) \leq \vol(B^d(r^* + \sqrt{d}/2))$. To prove the lower bound, we note that it is trivial for $r < \sqrt{d}/2$, and for $r \geq \sqrt{d}/2$, we can similarly set $r^* := r-\sqrt{d}/2$ and  $\vec{s}_2 := \vec{t}$ and then take $\vec{s}_1$ such that $N(d,(r^*)^2;\vec{s}_1) \geq \vol(B^d(r^*))$.
    
    It remains to prove \cref{eq:min_vs_average}. To that end, let
    \[
        f_{r'}(\vec{x}) := \begin{cases}
            1 &\|\vec{x}\| \leq r'\\
            0 &\text{otherwise}
        \end{cases}
    \]
    be the indicator function of $B^d(r')$.
    Then, we have
    \begin{align*}
        \vol(B^d(r')) 
            &= \int_{\R^d} f_{r'}(\vec{x}) \dl{\vec{x}} \\
            &= \sum_{\vec{z} \in \Z^d} \int_{[0,1)^d} f_{r'}(\vec{z} - \vec{s}) \dl{\vec{s}} \\
            &= \int_{[0,1)^d} \sum_{\vec{z} \in \Z^d} f_{r'}(\vec{z} - \vec{s}) \dl{\vec{s}} \\
            &= \int_{[0,1)^d} N(d,(r')^2;\vec{s}) \dl{\vec{s}}
        \; ,
    \end{align*}
    which immediately implies \cref{eq:min_vs_average}.
\end{proof}

\begin{remark} \cref{thm:volume_estimate} is simply a special case of a well known more general result, which says that for any full-rank lattice $\lat \subset \R^d$,
\[
    \vol(B^d(r - \mu(\lat)))/\det(\lat) \leq |\{\vec{y} \in \lat \ : \ \|\vec{y}\| \leq r\}| \leq \vol(B^d(r + \mu(\lat)))/\det(\lat)
    \; ,
\]
where $\mu(\lat) := \max_{\vec{t} \in \R^d} \min_{\vec{y} \in \lat} \|\vec{y} - \vec{t}\|$ is the lattice covering radius and $\det(\lat)$ is the lattice determinant. The proof of the more general result is essentially identical to the proof of \cref{thm:volume_estimate}, with the cube $[0,1)^d$ replaced by any fundamental domain of $\R^d/\lat$. This also can be extended to other norms by replacing the Euclidean ball with the unit ball of the norm in question.
\end{remark}

\begin{remark}
    Considerable effort has gone into studying the asymptotic behavior of the error term 
    \[
        h_d(r) := \vol(B^d(r)) - N(d,r^2) 
        \; 
    \]
    in \cref{eq:volume_approx} as the radius $r$ tends to infinity.

    In particular, Gauss's circle problem asks about this asymptotic behavior for $h_2$.
    Hardy and Landau showed that there exist arbitrarily large radii $r$ such that $|h_2(r)| \geq \widetilde{\Omega}(r^{1/2})$~\cite{hExpressionNumberSum1915,Lan15}, and Li and Yang showed that for all $r$ we have $|h_2(r)| \leq O(r^{\theta + \eps})$ for any $\eps > 0$ and $\theta \approx 0.629$~\cite{liImprovementGausssCircle2023}, slightly improving on earlier work of Huxley~\cite{huxleyExponentialSumsLattice2003a}. But the precise asymptotic behavior of $h_2(r)$ is still unknown. 
    The situation is similar for $h_3$. 
    
    For $d \geq 4$, the asymptotic behavior of the error $h_d$ is well understood, as it is known that $|h_d(r)| \leq \widetilde{O}(r^{d-2})$, and it is known that $|h_d(r)| \geq \widetilde{\Omega}(r^{d-2})$ for arbitrarily large radii $r$. 
    
    See \cite{ivicLatticePointsLarge2004} for a survey.
\end{remark}

\section{Bounds in terms of binomial coefficients, proven via induction}
\label{sec:binomials}

We begin our discussion by noting that the following simple recurrence holds. For any positive integers $d_1$, $d_2$, and $m$,
\begin{equation}
    \label{eq:r_d_recurrence}
    r_{d_1+d_2}(m)=\sum_{i=0}^{m}{r_{d_1}(i) \cdot r_{d_2}(m-i)}
    \; .
\end{equation}
This recurrence is identical to the binomial recurrence in \cref{lem:2dim_hockey_stick}.
This will in particular be used to prove the following theorem.

    \begin{theorem} 
    \label{thm:r_d_inequality_recurrence}
         For positive integers $d_1$ and $d_2$, $x_1,x_2 \in \R$, non-negative integers $y_1,y_2 \in \Z_{\geq 0}$ with $x_i + y_i \geq 0$, and $c_1,c_2, k > 0$, suppose that
        \[
       r_{d_i}(m) \leq  c_i k^m {m+x_i-1 \choose m-y_i}
        \; 
        \]
        holds simultaneously for all non-negative integers $m$.
        Then, for any non-negative integers $a_1$, $a_2$ with $a_1 + a_2 \geq 1$, and any non-negative integer $m$, we have
        \[
            r_{a_1d_1 + a_2 d_2}(m) \leq c_1^{a_1}c_2^{a_2} k^m {m+a_1x_1+a_2x_2-1 \choose m-a_1y_1-a_2y_2}
            \; .
        \]
        
        Similarly, if 
        \[
       r_{d_i}(m) \geq  c_i k^m {m+x_i-1 \choose m-y_i}
        \; 
        \]
        for all non-negative integers $m$, then
        \[
            r_{a_1 d_1 + a_2 d_2}(m) \geq c_1^{a_1}c_2^{a_2} k^m{m+a_1x_1+a_2 x_2-1 \choose m-a_1y_1-a_2 y_2}
            \; .
        \]
    \end{theorem}
    \begin{proof}
        Intuitively, the result follows from the fact that $r_{a_1 d_1 + a_2 d_2} (m)$ satisfies the same recurrence relation (i.e., \cref{eq:r_d_recurrence}) as the function $c_1^{a_1} c_2^{a_2} k^m \binom{m + a_1 x_1 + a_2 x_2 - 1}{ m -a_1y_1 - a_2 y_2}$.
        
        We will give a careful proof of the result about upper bounds. The analogous result for lower bounds has an essentially identical proof (after noting that since $x_i + y_i \geq 0$, all binomial coefficients encountered are non-negative). 
        
        We will prove that the  bound holds for all non-negative integers $a_1$ and $a_2$ via induction on $a_1+a_2$.
        For our base case, we will take the case where $a_1+a_2 = 1$, which is true by assumption.

        Now, consider an arbitrary such $a_1, a_2$ with $a_1 + a_2 > 1$. We may assume without loss of generality that $a_1 > 0$. We may therefore assume by induction that the result holds for $(a_1',a_2') := (a_1 - 1,a_2)$.
        Now, we can apply \cref{eq:r_d_recurrence} along with our induction hypothesis to get 
        \begin{align*}
            r_{a_1 d_1+a_2 d_2}(m) &= \sum_{i=0}^{m}{r_{d_1}(i)r_{(a_1-1)d_1+a_2 d_2}(m-i)}  \\
            &\leq \sum_{i=0}^{m} {c_1 k^i {i + x_1 - 1 \choose i-y_1} \cdot c_1^{(a_1-1)}c_2^{a_2}k^{m-i}{(m-i) + ((a_1-1)x_1+a_2x_2) -1 \choose (m-i)-(a_1-1)y_1-a_2y_2}} \\
            &= c_1^{a_1}c_2^{a_2} k^m \sum_{i=0}^{m}{{i + x_1 - 1 \choose i-y_1}{(m-i) + (a_1-1)x_1+a_2x_2 -1 \choose (m-i)-(a_1-1)y_1-a_2y_2}} \; .
        \end{align*} 
        
        Finally, by \cref{lem:2dim_hockey_stick}, we have that the final sum is equal to 
        \[ c_1^{a_1}c_2^{a_2} k^m {m + a_1x_1+a_2x_2 - 1 \choose m-a_1y_1-a_2y_2}
        \; ,
        \]
        as needed.
    \end{proof}

\subsection{A simple lower bound}

Combining our main theorem with the simple bounds for small $d$ in \cref{sec:simple_bounds_prelims} gives the following result, which is a slight improvement of the result in \cite{regevSimpleProofReverse2023}. 

\begin{theorem}
    \label{thm:simple_binomial_lower_bound}
    For all integers $d \geq 5$ and $m \geq 0$, we have 
    \[ 
        r_d(m) \geq {m + d/2 -1 \choose m}
        \; .
    \]
\end{theorem}
\begin{proof}
    By \cref{binomial_r_6_r_8_bounds}, the result holds for $d=6$ and $d=8.$ Applying \cref{thm:r_d_inequality_recurrence} gives us that our claim holds for even $d \geq 6$ such that $d \neq 10.$ 
    
    Using \cref{r_5_bound} and \cref{lem:r_7_and_r_9_lower_bound} we can apply \cref{thm:r_d_inequality_recurrence} three times with the cases $d_1=6$ and $d_2 =5,7,9$ where we choose $x_i = d_i/2$ and $y_i=0$ for $i \in \{1,2\}$ to obtain that the theorem holds for all odd integers $d$ such that $d \geq 5.$ 

    Finally, we can see that the theorem holds in the case of $d=10$ by applying \cref{thm:r_d_inequality_recurrence} where $d_1=d_2=5$ with $x_i = d_i/2$ and $y_i=0$ for $i \in \{1,2\}$ as before. 

    We have exhausted all positive integers $d \geq 5$ through casework, so the theorem is proven.
\end{proof}

\begin{corollary}
    \label{cor:simple_binomial_lower}
    For all integers $d \geq 5$ and $m \geq 0$, we have 
    \[ 
        N(d,m) \geq {m+d/2 \choose m}
        \; .
    \]
\end{corollary}
\begin{proof}
    We can apply \cref{thm:simple_binomial_lower_bound} and use \cref{lem:2dim_hockey_stick} in the case $x_1=d/2$ and  $x_2 = 1$ to get 
    \[ 
        N(d,m) = \sum_{i=0}^{m}{r_d(i)} \geq \sum_{i=0}^{m}{{i+d/2-1 \choose i}}  = \sum_{i=0}^{m}{{i+d/2-1 \choose i} \cdot {m-i + 1-1 \choose m-i}} = {m+d/2 \choose m}
        \; .
    \]
\end{proof}

\subsection{A simple universal upper bound}

\begin{theorem}
    \label{thm:simple_binomial_upper_bound}
    For all integers $d \geq 1$ and $m \geq 0$, we have
    \[ r_d(m) \leq \min \Big\{ {m+2d -1 \choose m }, 2^d \binom{m+d-1}{m}  \Big\}
    \; .\]
\end{theorem}
\begin{proof}
    Note that when $d=1$, we have that 
    \[ r_1(m) = \begin{cases}
        0 & m \text{ is not a square} \\
        1 & m = 0 \\
        2 & m > 0 \text{ is a square.}
    \end{cases}\]
    In particular, we have the simple upper bound $r_1(m) \leq m+1 = {m+1 \choose m}$.
    We can then apply \cref{thm:r_d_inequality_recurrence} in the case of $d_1=d_2=1$ with $x_1=x_2=2$ and $y_1=y_2=0$ to get 
    \[ 
        r_d(m) \leq {m+2d-1 \choose m}
    \; 
    ,\]
    as needed.

    Alternatively, we can apply the simple upper bound $r_1(m) \leq 2 = 2 \binom{m}{m}$ to obtain the bound $r_d(m) \leq 2^d \binom{m+d-1}{m}$. (As we noted in the introduction, the bound $r_d(m) \leq 2^d \binom{m+d-1}{m}$ also has a combinatorial proof. One can similarly obtain a bound of $r_d(m) \leq 2^m \binom{m+d-1}{m}$ via the same argument, but this bound is never better than the stated bound.)
\end{proof}

\begin{corollary}
    \label{cor:simple_binomial_upper}
    For all integers $d \geq 1$ and $m \geq 0$, we have  
    \[ 
        N(d,m) \leq \min\Big\{{m+2d \choose m}, 2^d \binom{m+d}{m} \Big\}
    \; .
    \]
\end{corollary}
\begin{proof}
    We can apply \cref{thm:simple_binomial_upper_bound} and \cref{lem:2dim_hockey_stick} with $x_1=2d$ and $x_2=1$ to obtain
    \[ 
        N(d,m) = \sum_{i=0}^{m}{r_d(i)} \leq \sum_{i=0}^{m}{{i+2d-1 \choose i}}= \sum_{i=0}^{m}{{i+2d-1 \choose i} \cdot {m-i+1-1 \choose m-i}} = {m+2d \choose m}
    \; ,
    \]
    as needed. The proof of the other bound is identical.
\end{proof}

\subsection{Some more complicated (and quite tight) bounds}

We now illustrate some of the more complicated bounds that one can obtain using \cref{thm:r_d_inequality_recurrence}. We stress that these are merely examples and that \cref{thm:r_d_inequality_recurrence} clearly implies an infinite family of bounds of this form (each obtained by starting with different bounds on $r_d(m) $ of the form $c k^m \binom{m+x-1}{m-y}$ for different choices of $d, c, k, x, y$). Nevertheless, the below bounds are quite tight, as illustrated in \cref{sec:some_plots}.

\begin{theorem}
    \label{thm:exponential_binomial}
    For all non-negative integers $a$ and $b$, let $d := 6a + 8b > 0$. Then, for any non-negative integer $m$, we have 
    \[ (20)^a(84)^b {m+a+2b-1 \choose m - 2a-2b} \leq r_d(m) \leq \min\{ 43^a 116^b, 1.04^m \cdot 21^a 44^b \} \cdot   {m + d/2 -1 \choose m} \; .\]
    (Note that the lower bound is zero when $m < 2a + 2b$, and that for fixed $d$ both bounds are strongest when $b$ is chosen to be maximal.)
\end{theorem}
\begin{proof}
    We claim that 
    \[
        20 \binom{m}{2} < r_6(m) < \min\{43, 21 \cdot 1.04^{m} \} \cdot \binom{m+2}{2}
        \; 
    \]
    for all non-negative integers $m$. For $m = 0$, $r_6(0) = 1$ and the bounds are immediate (noting that $\binom{0}{2} = 0$). For $m \geq 1$, the lower bound is immediate from \cref{r_6 bound}, as is the upper bound $r_6(m) < 21.5m^2 < 43\binom{m+2}{2}$. Furthermore, a quick check shows $43 < 21 \cdot 1.04^m$ for all $m \geq 19$. So, to finish proving the upper bound, it suffices to show that $r_6(m) < 21 \cdot 1.04^m \cdot \binom{m+2}{2}$ for all $1 \leq m \leq 18$, which is easy to check. (A slightly more careful bound gives $r_6(m) \leq (208/9) \cdot (15/13)^{(m-3)/4} \cdot \binom{m+2}{2} \approx 20.76 \cdot 1.036^m \cdot \binom{m+2}{2}$, which is tight for $m = 7$ and $m = 11$.)

    Similarly, we claim that
    \[
        84 \binom{m+1}{3} < r_8(m) < \min\{ 116, 44 \cdot 1.04^m\} \cdot \binom{m+3}{3}
    \]
    for all non-negative integers $m$.  Again, the lower bound and the upper bound of $116 \binom{m+3}{3} > 19.3m^3$ are immediate from \cref{r_8 bound}. To prove that $r_8(m) < 44 \cdot 1.04^m \cdot \binom{m+3}{3}$, it suffices to notice that this is weaker than the upper bound that we have already proven for $m \geq 25$. For $m \leq 24$, a quick direct check shows that the inequality holds.

    The upper and lower bounds on $r_d(m)$ then follow by applying \cref{thm:r_d_inequality_recurrence}. In particular, the lower bound on $r_d(m)$ follows by applying \cref{thm:r_d_inequality_recurrence} with $d_1 = 6$, $d_2 = 8$, $x_1 = 1$, $x_2 = 2$, $y_1 = y_2 = 2$, $k = 1$, $c_1 = 20$, and $c_2 = 84$. The upper bounds follow by taking $d_1 = 6$, $d_2 = 8$, $x_i = d_i/2$, $y_1 = y_2 = 0$, and $(c_1,c_2,k) \in \{(43,116,1), (21,44,1.04)\}$.
\end{proof}

The following corollary shows how to extend the lower bound above to odd integers $d \geq 5$.

\begin{corollary}
    For all non-negative integers $a, b, c$ and $d_0 \in \{5,7,9\}$ let $d := d_0 a + 6b + 8c > 0$. Then,
    \[ 
    r_d(m) \geq (20)^b(84)^c {m+d_0a/2+b+2c-1 \choose m - 2b-2c} 
    \; .\]
\end{corollary}
\begin{proof}
    By \cref{r_5_bound,lem:r_7_and_r_9_lower_bound}, we have
    \[ {m+d_0/2 -1 \choose m} \leq r_{d_0}(m)\]
    for $d_0 \in \{5,7,9\}$. Applying \cref{thm:r_d_inequality_recurrence} gives that the result holds when $b = c = 0$. So, we assume that $b + c \geq 1$ in what follows.

    By  \cref{thm:exponential_binomial}, we know that 
    \[ (20)^b(84)^c {m+b+2c-1 \choose m - 2b-2c} \leq r_{6b+8c}(m) \;.\]
    Then, we obtain the result by applying \cref{thm:r_d_inequality_recurrence} in the case of $d_1 = 6b+8c$ and $d_2 = d_0$ with $x_1 = b+2c, y_1 = 2b+2c$ and $x_2 = d_0/2$, $y_2 = 0$, $c_1 =20^b84^c $, and $c_2 = 1$. 
\end{proof}

\section{Bounds based on the theta function in the style of Mazo and Odlyzko}

Recall that for $\tau > 0$ and $\vec{t} = (t_1,\ldots, t_d) \in \R^d$, we write 
\[
    \Theta(\tau;\vec{t}) := \sum_{\vec{z} \in \Z^d} q^{\|\vec{z} - \vec{t}\|^2} = \prod_{i=1}^d \Theta(\tau;t_i)
    \; ,
\]
where $q := e^{-\pi \tau}$. In this section, we prove the following result.

\begin{theorem}
    \label{thm:main_theta_theorem}
    For any $\vec{t} \in \R^d$ and $m \geq 1$, we have
    \[
         N(d,m;\vec{t}) = e^{-\chi(m;\vec{t})}\inf_{\tau > 0} \frac{\Theta(\tau;\vec{t})}{q^{m}}
         \; ,
    \]
    where $q := e^{-\pi \tau}$ and $\chi(m;\vec{t})$ satisfies
    \[
        0 \leq \chi(m; \vec{t}) < 3+\min \{4\sqrt{d}, 8 \sqrt{m}\log'(d/m)\}
        \; ,
    \]
    where $\log' x := \max\{1,\log x\}$.
\end{theorem}

To prove \cref{thm:main_theta_theorem}, we work with $D(\tau;\vec{t})$. Recall that $D(\tau;\vec{t})$ is the probability distribution over $\Z^d - \vec{t}$ induced by the Gaussian measure, known as the discrete Gaussian distribution. In other words, for any $\vec{x} \in \Z^d - \vec{t}$, we have
\[
    \Pr_{\vec{X} \sim D(\tau; \vec{t})}[\vec{X} = \vec{x}] := \frac{q^{\|\vec{x}\|^2}}{\Theta(\tau;\vec{t})}
    \; .
\]
The following basic lemma shows that, in order to prove tight bounds on $N(d,m;\vec{t})$, it suffices to prove that for some $\tau > 0$, $D(\tau;\vec{t})$ is concentrated on vectors with squared norm roughly $m$. We will then prove the desired concentration by carefully bounding the variance of $\|\vec{X}\|^2$ for $\vec{X} \sim D(\tau;\vec{t})$ (separately for $\tau \geq 1$ and $\tau \leq 1$) and simply applying Chebyshev's inequality.

\begin{lemma}
    \label{lem:markov_theta}
    For any $m \geq 0$, $\vec{t} \in \R^d$, $\tau > 0$, and $\delta \geq 0$, 
    \begin{equation}
        \label{eq:Markov_theta_and_reverse_Markov_theta}
        \Pr_{\vec{X} \sim D(\tau;\vec{t})}[m-\delta \leq \|\vec{X}\|^2\leq m] \cdot \frac{\Theta(\tau;\vec{t})}{q^{m-\delta}} \leq N(d,m;\vec{t})\leq \frac{\Theta(\tau;\vec{t})}{q^{m}}
        \; .
    \end{equation} 
\end{lemma}
\begin{proof}
     The upper bound simply follows from the fact that $q^{\|\vec{x}\|^2} \geq q^m$ when $\|\vec{x}\|^2 \leq m$. 
   To prove the lower bound, notice that it suffices to consider the case when $0 \leq \delta \leq m$ (since the $\delta > m$ case follows trivially from the $\delta = m$ case). Then, we have
      \begin{align*}
       N(d,m;\vec{t})
       &\geq\frac{1}{q^{m-\delta}}\sum_{\substack{m-\delta
       \leq \|\vec{z}-\vec{t}\|^2
       \leq m
    \\\vec{z}\in\Z^d }}q^{\|\vec{z}-\vec{t}\|^2}\\
    &=\Pr_{\vec{X} \sim D(\tau;\vec{t})}[m-\delta \leq \|\vec{X}\|^2\leq m] \cdot \frac{\Theta(\tau;\vec{t})}{q^{m-\delta}}
    \; ,
   \end{align*}
   as needed.
\end{proof}

        It is not immediately clear how tight \cref{lem:markov_theta} is because it is not immediately clear how to work with the probability in the lower bound. We would therefore like to replace the lower bound in \cref{lem:markov_theta} with a bound that is easier to work with, and ideally one that is evidently close to the upper bound. Whether this is possible of course strongly depends on the choice of $\tau$ and $\delta$. Indeed, for choosing $\tau$, it is quite natural to consider the special case when $m \approx \expect_{\vec{X} \sim D_{(\tau;\vec{t})}}[\|\vec{X}\|^2]$, as we do below.\footnote{Indeed, it is not hard to see that  $\Theta(\tau;\vec{t})/q^m$ is minimized at $\tau$ if and only if $\tau$ satisfies $m = \expect_{\vec{X} \sim D_{(\tau;\vec{t})}}[\|\vec{X}\|^2]$.} And, it is natural to take $\delta$ to be roughly equal to the standard deviation of $\|\vec{X}\|^2$ when $\vec{X} \sim D(\tau;\vec{t})$. Then, we can simply apply Chebyshev's inequality to bound the probability in the lower bound in \cref{lem:markov_theta}. The below makes this precise. (Notice that here we think of $\tau$ fixed and $m$ chosen depending on $\tau$, but we will quickly switch to fixing $m$ and choosing $\tau$ depending on $m$.)

\begin{lemma}
    \label{thm: application of Chebyshev's Inequality}
    For any $\vec{t} \in \R^d$, $\tau > 0$, and $c > 1$, we have
    \[\Pr_{\vec{X}\sim D(\tau;\vec{t})}[m-\delta \leq \|\vec{X}\|^2\leq m]\geq 1-\frac{1}{c^2} 
    \; ,
    \]
    where 
    \[\delta := \delta(\tau;\vec{t})
        :=2c\sqrt{\Var_{\vec{X}\sim D(\tau;\vec{t})}[\|\vec{X}\|^2]}\; ,\]
    and 
    \[ m := m(\tau; \vec{t})
    :=\expect_{\vec{X}\sim D(\tau;\vec{t})}[\|\vec{X}\|^2]+\frac{\delta}{2}\;.\]
\end{lemma}

\begin{proof}
    We apply Chebyshev's inequality and get
    \[
    \Pr_{\vec{X}\sim D(\tau;\vec{t})}\biggl[\biggl|\|\vec{X}\|^2-\Big(m-\frac{\delta}{2}\Big)\biggr|\geq \frac{\delta}{2}\biggr]=\Pr_{\vec{X}\sim D(\tau;\vec{t})}\biggl[\bigl|\|\vec{X}\|^2-\expect[\|\vec{X}\|^2]\bigr|\geq \frac{\delta}{2}\biggr]
    \leq \frac{\Var[\|\vec{X}\|^2]}{(\delta/2)^2}= \frac{1}{c^2}\;.\]
    So,
    \[
    \Pr_{\vec{X}\sim D(\tau;\vec{t})}[m-\delta \leq \|\vec{X}\|^2\leq m]\geq 1-\frac{1}{c^2} \;,
    \]
    as needed.
\end{proof}

\begin{corollary}
    \label{cor:theta_infimum}
    For any $\vec{t} \in \R^d$, $m > 0$, and $c > 1$, we have
    \[
        (1-1/c^2) \cdot \inf_{\tau > \tau^*} \frac{\Theta(\tau;\vec{t})}{q^{m-c \eps}} \leq N(d,m;\vec{t}) \leq \inf_{\tau > 0} \frac{\Theta(\tau;\vec{t})}{q^m}
        \; ,
    \]
    where  $q := e^{-\pi \tau}$, $\tau^* := 0$ if $d \leq 2m$ and $\tau^* := \log(2(d-m)/m)/\pi$ if $d > 2m$,\footnote{The fact that we can restrict the infimum in the lower bound to $\tau > \tau^*$ rather than allowing any $\tau > 0$ should be viewed as a rather minor optimization. If we allowed arbitrary $\tau > 0$ here, our bound on $\chi$ in \cref{thm:main_theta_theorem} would go from roughly $\sqrt{m} \log(d/m)$ to roughly $\sqrt{m} \log(d/m) + \log^2(d/m)$. By restricting the infimum to $\tau > \tau^*$ in the lower bound, we are able to avoid this issue and to simplify the proof of \cref{thm:main_theta_theorem}.} and
    \[\eps := \eps(\tau;\vec{t})
        := 
        2\sqrt{\Var_{\vec{X}\sim D(\tau;\vec{t})}[\|\vec{X}\|^2]}\; .
    \]
    
    Furthermore, if $m \leq \min_{\vec{z} \in \Z^d} \|\vec{z} - \vec{t}\|^2$, then 
    \[
        N(d,m;\vec{t}) = \inf_{\tau > 0} \frac{\Theta(\tau;\vec{t})}{q^m} = \lim_{\tau \to \infty} \frac{\Theta(\tau;\vec{t})}{q^m}
        \; .
    \]
\end{corollary}
\begin{proof}
    The upper bound is immediate by taking an infimum over all $\tau > 0$ of the upper bound in \cref{lem:markov_theta}. 
    
    For the lower bound, let $\delta := \delta(\tau;\vec{t}) := c \eps$ and $E(\tau;\vec{t}) := \expect_{\vec{X} \sim D(\tau;\vec{t})}[\|\vec{X}\|^2]$. We consider two cases. First, suppose that there exists a $\tau > 0$ such that $m = E(\tau;\vec{t}) + \delta(\tau;\vec{t})/2$. Then, for this particular value of $\tau$, we can apply \cref{thm: application of Chebyshev's Inequality} to conclude that
    \[
        \Pr_{\vec{X}\sim D(\tau;\vec{t})}[m-\delta \leq \|\vec{X}\|^2\leq m]\geq 1-\frac{1}{c^2} 
        \; .
    \]
	Plugging this in to \cref{lem:markov_theta}, we see that for this particular choice of $\tau$,
    \[
        N(d,m;\vec{t}) \geq (1-1/c^2) \frac{\Theta(\tau;\vec{t})}{q^{m-\delta}}
        \; .
    \]
    
    The lower bound follows in this case, provided that we can show that such a $\tau$ must satisfy $\tau > \tau^*$. This is trivial when $d \leq 2m$ (since then $\tau^* = 0$), so we assume $d > 2m$. Indeed, by \cref{thm:second_moment_lower_bound}, we have that $E(\tau;\vec{t}) \geq 2dq/(1+2q)$.
    In particular, if $\tau \leq \tau^*$, we see that $E(\tau;\vec{t}) \geq m$, contradicting the fact that $E(\tau;\vec{t}) + \delta(\tau;\vec{t})/2 = m$. So, we must have $\tau > \tau^*$, as needed.

    Now, suppose that no such $\tau$ exists.\footnote{In some sense, this case is not very interesting, since as we will see, this case simply corresponds to the case when $m \leq \min_{\vec{z} \in \Z^d} \|\vec{z} - \vec{t}\|^2$, in which case $N(d,m;\vec{t})$ is trivial to compute.} For this, notice that the functions $\delta(\tau;\vec{t})$ and $E(\tau;\vec{t})$ are continuous (in fact, analytic) functions of $\tau > 0$. Therefore, the function $f(\tau) := E(\tau;\vec{t}) + \delta(\tau;\vec{t})/2$ is itself continuous. We also have that $\lim_{\tau \to 0}f(\tau) = \infty$ and 
    \[
        \lim_{\tau \to \infty} f(\tau) = \min_{\vec{z} \in \Z^d} \|\vec{z}- \vec{t}\|^2
        \; .
    \]
    It follows that the image of $f$ over the positive real numbers includes all numbers that are strictly greater than $m^* := \min_{\vec{z} \in \Z^d} \|\vec{z} - \vec{t}\|^2$. Since by assumption no $\tau > 0$ exists such that $f(\tau) = m$, we must have $m \leq m^*$. So, we have reduced to the ``furthermore'' in this case.

	So, it remains to show the ``furthermore.''
To that end, notice that we have 
\[\lim_{\tau \to \infty }\Theta(\tau;\vec{t})/q^{m^*}  = N(d, m^*;\vec{t})
\; .
\]
Therefore, we have
    \[
        \lim_{\tau \to \infty} \frac{\Theta(\tau;\vec{t})}{q^{m}}  = N(d,m^*;\vec{t}) \cdot \lim_{q \to 0}q^{m^*-m}
        \; .
    \]
    So, if $m = m^*$, we clearly have that the limit is equal to $N(d,m;\vec{t})$ as needed, and if $m < m^*$, then the limit is equal to zero and $N(d,m;\vec{t})$ is also equal to zero. In either case, we see that 
    \[
     \lim_{\tau \to \infty} \frac{\Theta(\tau;\vec{t})}{q^{m}} = N(d,m;\vec{t})
     \; .
    \]
     The fact that the infimum is equal to the limit in this case then follows from the fact that
    \[
        \inf_{\tau > 0} \frac{\Theta(\tau;\vec{t})}{q^{m}} \leq \lim_{\tau \to \infty} \frac{\Theta(\tau;\vec{t})}{q^{m}} = N(d,m;\vec{t}) \leq \inf_{\tau > 0} \frac{\Theta(\tau;\vec{t})}{q^{m}}
        \; ,
    \]
    where the last inequality is the upper bound shown earlier.
    \end{proof}

Of course, it is not immediately clear from \cref{cor:theta_infimum} how tight these bounds are or whether they are useful. The rest of this section is devoted to showing that they are relatively tight. In particular, by bounding the variance of $D(\tau;\vec{t})$, we will prove \cref{thm:main_theta_theorem}. (Along the way, we observe that the upper and lower bounds in \cref{cor:theta_infimum} can be computed efficiently to high accuracy by simply truncating either the defining series or its dual. Of course, this is well known and $\Theta$ is implemented in many computer algebra systems.)

\subsection{Approximating \texorpdfstring{$\Theta$}{the theta function} and its associated moments by finite sums}
\label{sec:finite_sum_defs}

In this section, we observe that the function $\Theta(\tau;\vec{t})$ and the associated moments $\expect_{\vec{X} \sim D(\tau;\vec{t})}[\|\vec{X}\|^2]$ and $\Var_{\vec{X} \sim D(\tau;\vec{t})}[\|\vec{X}\|^2]$ can be expressed as rapidly convergent series and bound the error obtained by truncating such series. This in particular means that the bounds in \cref{cor:theta_infimum} can be computed to arbitrary accuracy. (Of course, many computer algebra packages have these functions implemented.) This will also be useful because it will allow us to obtain a simple bound on the error in \cref{cor:theta_infimum}. (In fact, we will only need the simple $k = 1$ case to prove our main theorem, but we prove a more general---though still quite straightforward---result because we expect it to be useful.)

Since $\Theta(\tau;\vec{t}) = \prod_{i=1}^{d} \Theta(\tau;t_i)$ and since (more-or-less equivalently) $D(\tau;\vec{t})$ is a product distribution, in order to approximate $\Theta(\tau;\vec{t})$ and the moments of $D(\tau;\vec{t})$, it suffices to show how to approximate these in the one-dimensional case when $d = 1$.

If, say, $\tau \geq 1$, then this is more-or-less immediate. In particular, for any $\tau > 0$, $t \in \R$, and integer $k \geq 1$, we define
\[
    \Theta_{k}(\tau;t) :=  \sum_{z=1}^{k} q^{(z-t)^2} + \sum_{z=0}^{k-1} q^{(z+t)^2}
        \; , \quad \qquad \mu_k(\tau;t) := \sum_{z=1}^{k} (z-t)^2 q^{(z-t)^2} + \sum_{z=0}^{k-1} (z+t)^2 q^{(z+t)^2}
        \; ,
\]
and
\[
M_k(\tau;t) := \sum_{z=1}^{k} (z-t)^4 q^{(z-t)^2} + \sum_{z=0}^{k-1} (z+t)^4 q^{(z+t)^2}
\; ,
\]
where as always, $q := e^{-\pi \tau}$. We also define \[
\mu(\tau;t) := \lim_{k \to \infty} \mu_k(\tau;t)
\; , \text{ and } \quad M(\tau;t) := \lim_{k \to \infty} M_k(\tau;t)
\]
(just like $\Theta(\tau;t) = \lim_{k \to \infty} \Theta_k(\tau;t)$).

The following lemma shows that these finite sums yield good approximations to the relevant series when $k$ is not too small relative to $1/\sqrt{\tau}$. (These bounds will also be useful in the next section.)

As the proof is rather tedious and unilluminating, we relegate it to \cref{app:finite_sums_primal}. (Note that we did not put much effort into optimizing the constants, since the error terms are already quite small even without such effort.)

\begin{lemma}
    \label{lem:finite_sums_primal_new}
    For any $\tau > 0$, $t \geq 0$, and integer $k \geq  t + \sqrt{2/(\pi \tau)}$, we have
    \[
        0 \leq \Theta(\tau;t) - \Theta_k(\tau;t) \leq q^{k^2} + \frac{q^{w^2}}{\pi \tau w}
        \; ,
    \]
     \[
     0 \leq \mu(\tau;t) -  \mu_k(\tau;t) \leq k^2 q^{k^2} +  \frac{5w}{4\pi \tau} \cdot q^{w^2}
        \; ,
    \]
    and
    \[
        0 \leq  M(\tau;t) -  M_k(\tau;t) \leq k^4 q^{k^2} + \frac{31 w^3}{16 \pi \tau} q^{w^2} 
        \; ,
    \]
    where $w := k-t$ and $q := e^{-\pi \tau}$.
\end{lemma}

From this, we derive the following convenient corollary, which focuses on the case $k = 1$ and $t \in[0,1/2]$. This bound  will be useful for bounding the variance of $D(\tau;\vec{t})$. The proof simply applies \cref{lem:finite_sums_primal_new} with $k = 2$ and bounds the difference between $\Theta_2$ and $\Theta_1$ (and similarly for $\mu$ and $M$). We again relegate the proof to \cref{app:finite_sums_primal}.

\begin{corollary}
    \label{cor:k=1}
     For any $\tau \geq 1$ and $t \in [0,1/2]$, we have
     \[
         0 \leq \Theta(\tau;t) - \Theta_1(\tau;t) < 1.1q^{1+t^2}
         \; ,
     \]
      \[
      0 \leq \mu(\tau;t) -  \mu_1(\tau;t) < 1.1 q^{1+t^2}
         \; ,
     \]
     and
     \[
         0 \leq  M(\tau;t) -  M_1(\tau;t) < 1.1 q^{1+t^2}
         \; ,
     \]
    where $q := e^{-\pi \tau}$.
\end{corollary}

When $\tau$ is very small, $\Theta_k$, $\mu_k$, and $M_k$ are not actually good approximations to $\Theta$, $\mu$, and $M$ (except for very large $k \gtrsim 1/\sqrt{\tau}$), because the respective series do not converge rapidly when $\tau$ is small. In \cref{app:finite_sums_dual}, we show similar approximations that work when $\tau$ is small using the Poisson summation formula.

    \subsection{Bounds on the variance of \texorpdfstring{$D(\tau;\vec{t})$}{the discrete Gaussian}}

    It remains to show that the bound that we have proved on $N(d,m;\vec{t})$ in \cref{cor:theta_infimum} is relatively tight. 

    Notice that the upper and lower bounds in \cref{cor:theta_infimum} are the same except for the appearance of $\eps$ in the lower bound (and the rather minor effect of the $1-1/c^2$ term). So, to show that the bounds are quite tight, it suffices to bound $\eps := 2 \sqrt{\Var[\|\vec{X}\|^2]}$. We separately bound the variance in the case when $\tau$ is small and when $\tau$ is large.

        For $\tau \leq 1$, we obtain a quite tight bound on the variance by applying the Poisson summation formula, as below. (Note that $d/(2\pi^2 \tau^2)$ corresponds to the variance of the continuous Gaussian distribution, and it therefore arises naturally in this context.)

    \begin{lemma}
    \label{lem:var_bound_small_tau_dual}
    For any $\vec{t} \in \R^d$ and $0 < \tau \leq 1$, we have
    \[
        \Var_{\vec{X} \sim D(\tau;\vec{t})} [\|\vec{X}\|^2] \leq \frac{d}{2\pi^2 \tau^2} + \frac{5 d p}{\tau^{4}}
        \; ,
    \]
    where $p := e^{-\pi/\tau}$.
\end{lemma}
\begin{proof}
    Let
        \[
            V(\tau;t) := \Var_{X \sim D(\tau;t)}[X^2] = \frac{M(\tau;t)}{\Theta(\tau;t)} - \frac{\mu(\tau;t)^2}{\Theta(\tau;t)^2}
        \]
        be the variance of the one-dimensional distribution $D(\tau;t)$, where $M(\tau;t)$ and $\Theta(\tau;t)$ are as defined in \cref{sec:finite_sum_defs}.
        Since the $d$-dimensional distribution $D(\tau;\vec{t})$ is a product distribution, it follows that
        \[
            \Var_{\vec{X} \sim D(\tau;\vec{t})}[\|\vec{X}\|^2] = \sum_{i=1}^d V(\tau;t_i)
            \; .
        \]
        So, it suffices to prove that 
        \[
            V(\tau;t) \leq \frac{1}{2\pi^2 \tau^2} + \frac{5 p}{\tau^{4}}
            \; ,
        \]
    for all $t \in \R$.
    
    By \cref{fact:PSF_theta_mu_M} (i.e., the Poisson summation formula applied to $\Theta$, $\mu$, and $M$), we have
    \begin{align*}
       V(\tau;t) 
        &= \frac{1}{2 \pi^2 \tau^2} - \frac{4}{\pi \tau^{7/2}} \cdot \frac{S_2}{\Theta(\tau;t)} + \frac{2}{\tau^{9/2}} \cdot \frac{S_4}{\Theta(\tau;t)} - \frac{4}{\tau^{5}} \cdot \frac{S_2^2}{\Theta(\tau;t)^2}\\
        &\leq \frac{1}{2 \pi^2 \tau^2} + \frac{4}{\pi \tau^{7/2}} \cdot \frac{|S_2|}{\Theta(\tau;t)} + \frac{2}{\tau^{9/2}} \cdot \frac{|S_4|}{\Theta(\tau;t)}
       \; ,
    \end{align*}
    where we write
    \[
        S_2 := \sum_{z=1}^\infty \cos(2\pi zt) z^2 p^{z^2} 
        \; ,
    \]
    and
    \[
        S_4 := \sum_{z=1}^\infty \cos(2\pi zt) z^4 p^{z^2}
        \; .
    \]
   In particular,
    \[
        |S_2| \leq \sum_{z=1}^\infty z^2 p^{z^2} = p + \sum_{z=2}^\infty z^2 p^{z^2} \leq \Big(1+\frac{\tau}{2\pi} + \frac{\tau^2}{4\pi^2 }\Big) p < 1.2 p
        \; ,
    \]
    where the second-to-last inequality is \cref{clm:Gaussian_sum_bounds} and the last inequality uses the fact that $\tau \leq 1$. Similarly,
    \[
        |S_4| \leq \sum_{z=1}^\infty z^4 p^{z^2} = p + \sum_{z=2}^\infty z^4 p^{z^2} \leq \Big(1 + \frac{\tau}{2\pi} + \frac{3\tau^2}{4\pi^2} + \frac{3\tau^3}{8\pi^3} \Big) \cdot p < 1.25p
        \; .
    \]

     Finally, notice that by applying \cref{fact:PSF_theta_mu_M} again, we have
     \[
        \Theta(\tau;t) = \frac{1}{\sqrt{\tau}} + \frac{2}{\sqrt{\tau}} \cdot \sum_{z=1}^\infty \cos(2\pi t z)p^{z^2} \geq \frac{1}{\sqrt{\tau}} - \frac{2}{\sqrt{\tau}} \cdot \sum_{z=1}^\infty p^{z^2}  \geq \frac{1}{\sqrt{\tau}} - \frac{2p}{\sqrt{\tau}} - \frac{p \sqrt{\tau}}{\pi} > \frac{0.89}{\sqrt{\tau}}
        \; ,
     \]
     where the second-to-last inequality is again \cref{clm:Gaussian_sum_bounds} and the last inequality again uses the fact that $\tau \leq 1$.
     
     We conclude that
     \[
        \frac{M(\tau;t)}{\Theta(\tau;t)} - \frac{\mu(\tau;t)^2}{\Theta(\tau;t)^2} < \frac{1}{2 \pi^2 \tau^2} + \Big( 1.20 \cdot \frac{4}{\pi \tau^{7/2}} + 1.25 \cdot \frac{2}{\tau^{9/2}} \Big) \cdot \frac{\sqrt{\tau} p}{0.89} < \frac{1}{2\pi^2 \tau^2} + \frac{5p}{\tau^{4}}
        \; ,
     \]
     as claimed.
\end{proof}

We now handle the case when $\tau \geq 1$ by appropriately truncating the summations defining the variance. (The expression $4d/(e^2 \pi^2 \tau^2)$ arises naturally as a bound on the variance of the distribution obtained by restricting $D(\tau;\vec{t})$ to $\{0,1\}^d - \vec{t}$, and it is tight up to a factor of $2$. The second bound of the form $4\|\vec{t}\|^2/\tau^2 + 6dq$ is not particularly natural or tight in general. However, the bound is only relevant for us when $\|\vec{t}\|^2 \lesssim dq \tau^2$. In that case, the bound simply says that the variance is $O(dq)$, which is in fact tight.)

    \begin{lemma}
        \label{lem:var_bound_big_tau_primal}
        For any $\vec{t} \in \R^d$ and $\tau \geq 1$, we have
        \[
            \Var_{\vec{X} \sim D(\tau;\vec{t})}[\|\vec{X}\|^2] < \min\Big\{\frac{4d}{e^2 \pi^2 \tau^2} + 3dq, \frac{4\|\vec{t}\|^2}{\tau^2} + 6dq\Big\} 
            \; ,
        \]
        where $q := e^{-\pi \tau}$.
    \end{lemma}
    \begin{proof}
    We may assume that $\vec{t} \in [0,1/2]^d$ without loss of generality.

        Let
        \[
            V(\tau;t) := \Var_{X \sim D(\tau;t)}[X^2]
        \]
        be the variance of the one-dimensional distribution $D(\tau;t)$.
        Since the $d$-dimensional distribution $D(\tau;\vec{t})$ is a product distribution, it follows that
        \[
            \Var_{\vec{X} \sim D(\tau;\vec{t})}[\|\vec{X}\|^2] = \sum_{i=1}^d V(\tau;t_i)
            \; .
        \]
        So, it suffices to prove that
        \[
            V(\tau;t) < \min \Big\{ \frac{4}{e^2 \pi^2 \tau^2} + 3q, \frac{4t^2}{\tau^2} + 6q  \Big\}
        \]
        for $t \in [0,1/2]$.

        To that end, using the notation from \cref{sec:finite_sum_defs}, we have
        \[
            V(\tau;t) = \frac{1}{\Theta(\tau;t)^2} \cdot \big( \Theta(\tau;t)M(\tau;t)  - \mu(\tau;t)^2\big)
            \; .
        \]
        Applying \cref{cor:k=1} and the trivial fact that $\Theta(\tau;t)^2 \geq q^{2t^2}$, we see that
        \begin{align*}
            V(\tau;t) 
                &\leq \frac{1}{q^{2t^2}} \cdot \big(  ( \Theta_1(\tau;t) + 1.1q^{1+t^2})  (M_1(\tau;t) + 1.1 q^{1+t^2} ) - \mu_1(\tau;t)^2 \big) 
            \; .
        \end{align*}
        Plugging in the expressions for $\Theta_1$, $\mu_1$, and $M_1$ and simplifying, we see that
        \[
            V(\tau;t) \leq  q^{1-2t} (1-2t)^2 + 1.1 q^{1-t^2} \cdot M_1(\tau;t) + 1.1q^{1-t^2} \Theta_1(\tau;t)  + 1.21 q^{2}
            \; .
        \]
        Then, using the fact that $\Theta_1(\tau;t) = q^{t^2} + q^{(1-t)^2} \leq 2 q^{t^2}$ and the fact that $M_1(\tau;t) = t^4 q^{t^2} + (1-t)^4 q^{(1-t)^2} \leq q^{t^2}/8$,\footnote{The fact that $M_1(\tau;t) \leq q^{t^2}/8$ can be seen by noticing that $f(t) := M_1(\tau;t)/q^{t^2} = t^4 + (1-t)^4 q^{1-2t}$ is an increasing function of $t$ for $0 \leq t \leq 1/2$. Indeed, $f'(t) = 4t^3 + 2(1-t)^3 q^{1-2t}(\pi \tau (1-t) -2)$. If $(\pi \tau (1-t) -2) > 0$, then clearly $f'(t) > 0$, and we are done. Otherwise, we have $f'(t) \geq 4t^3 + 2(1-t)^3 (\pi (1-t) -2)$, and a quick check shows that this polynomial in $t$ has no real roots and therefore is always positive, as needed.} we see that
        \[
            V(\tau;t) \leq q^{1-2t} (1-2t)^2 + (1.1/8 + 2.2+ 1.21 q) \cdot q < q^{1-2t} (1-2t)^2 + 3q
            \; .
        \]

        The bound $V(\tau;t) < \frac{4}{e^2 \pi^2 \tau^2} + 3q$ then follows immediately from the fact that $q^{1-2t} (1-2t)^2$ achieves its unique maximum in the interval $t \in [0,1/2]$ at the point $t = 1/2 - 1/(\pi \tau)$, where it takes the value $4/(e^2 \pi^2 \tau^2)$. 
        
        To prove the bound $V(\tau;t) < 4t^2/\tau^2 + 6q$, first, notice that if $t \geq t_1 := 1/4$, then this bound is strictly weaker than the bound that we have already proved.  Next, notice that if $t \leq t_0 := 1/(2\pi \tau)$, then we see that $V(\tau;t) \leq q^{1-2t_0} + 3q = eq + 3q < 6q$, as needed. 
        So, in what follows, we may assume $t_0 < t < t_1$.
       Since in particular this means that $t > 0$, we may define
        \[g(t) := q^{1-2t} (1-2t)^2/t^2
        \; ,
        \]
        and
        we see that $V(\tau;t) \leq g(t) \cdot t^2 + 3q$. So, it suffices to prove that 
        \[g(t) \leq \frac{4}{\tau^2}
        \]
        for $t \in [t_0,t_1]$. 
        
        By differentiating $g$, we see that $g$ has no local maxima on the interval $t \in [t_0,t_1]$.\footnote{In more detail, if $\tau \leq 8/\pi$, then $g(t)$ is a strictly decreasing function on the entire interval $t \in (0,1/2]$. If $\tau > 8/\pi$, then $g(t)$ is increasing in the interval $t \in [(1-\eta)/4, (1+\eta)/4]$, where $\eta := \sqrt{1-8/(\pi \tau)}$. However, notice that $(1+\eta)/4 > 1/4 = t_1$. (This is why we chose this value of $t_1$.)} So, we conclude that $g(t) \leq \max\{g(t_0), g(t_1)\}$. Indeed, we have $g(t_0) = 4eq(\pi \tau - 1)^2 < 4/\tau^2$.\footnote{The inequality here follows from the fact that the ratio of the left-hand side to the right-hand side is maximized for $\tau \geq 1$ at $\tau = (5 + \sqrt{17})/(2\pi)$, where it is equal to roughly $0.76 < 1$.} And, 
        \[
            g(t_1) = (1-2t_1)^2 q^{1-2t_1}/t_1^2 = 4\sqrt{q} < 4/\tau^2
            \; ,
        \]
        as needed.
    \end{proof}

\subsection{Finishing the proof\texorpdfstring{ of \cref{thm:main_theta_theorem}}{}}

Finally, we plug our bounds on the variance into \cref{cor:theta_infimum} to obtain \cref{thm:main_theta_theorem}.

\begin{proof}[Proof of \cref{thm:main_theta_theorem}]
    The fact that $\chi(m;\vec{t}) \geq 0$ is immediate from the upper bound in \cref{cor:theta_infimum}. And, by the ``furthermore'' of \cref{cor:theta_infimum}, we can assume that $m > \min_{\vec{z} \in \Z^d}\|\vec{z} - \vec{t}\|^2$.
    
    For the upper bound on $\chi$, let 
    \[
        f(\tau;\vec{t}) := \log \Big( \frac{\Theta(\tau;\vec{t})}{q^{m}}\Big)
        \; .
    \]
    By \cref{cor:theta_infimum}, we have that
    \begin{equation*}
        \label{eq:lower_bound_inf_in_proof}
        N(d,m;\vec{t}) \geq (1-1/c^2)\cdot \inf_{\tau > \tau^*} q^{c \eps} e^{f(\tau;\vec{t})}
        \; ,
    \end{equation*}
    where, say, $c := 30/29$, 
    \[
        \eps := \eps(\tau;\vec{t}) := 2\sqrt{\Var_{\vec{X}\sim D(\tau;\vec{t})}[\|\vec{X}\|^2]}
        \; ,
    \]
    and $\tau^* := 0$ if $d \leq 2m$  and otherwise $\tau^* := \log(2(d-m)/m)/\pi$.
    So, we have
    \begin{align}
        \chi(m;\vec{t}) 
            &\leq -\log(1-1/c^2) +  \inf_{\tau > \tau^*} ( f(\tau; \vec{t}))  - \inf_{\tau >\tau^*} \big(f(\tau;\vec{t}) - c \pi \tau \eps(\tau;\vec{t}) \big) \nonumber \\
            &< 3+ c \pi \sup_{\tau > \tau^*} \tau \eps(\tau;\vec{t}) \label{eq:sup_of_eps}
        \; .
    \end{align}
    It remains to bound this supremum.

    We first prove that $\chi < 3 + 4\sqrt{d}$ by showing that $c \pi \tau \eps(\tau;\vec{t}) \leq 4\sqrt{d}$ for all $\tau > 0$ and all $\vec{t} \in \R^d$.
    Indeed, if $\tau \leq 1$, then it follows from \cref{lem:var_bound_small_tau_dual} that 
    \[
        \pi \eps(\tau;\vec{t}) \leq \frac{\sqrt{2d}}{\tau} \cdot \sqrt{ 1 + \frac{10\pi^2 p}{\tau^2}} \leq \frac{\sqrt{2d}}{\tau} \cdot  \sqrt{1+10\pi^2 e^{-\pi}} < \frac{3.3\sqrt{d}}{\tau }
        \; .
    \]
    On the other hand, if $\tau \geq 1$, then it follows from \cref{lem:var_bound_big_tau_primal} that
    \[
        \pi \eps(\tau;\vec{t}) \leq \frac{4 \sqrt{d}}{e \tau} \cdot \sqrt{1+3e^2 \pi^2 \tau^2 q/4} \leq \frac{4 \sqrt{d}}{e \tau} \cdot \sqrt{1+3e^2 \pi^2 e^{-\pi}/4} < \frac{3 \sqrt{d}}{ \tau}
        \; .
    \]
    Either way, we have $\pi c \tau \eps(\tau;\vec{t}) < 3.3c \sqrt{d} < 4 \sqrt{d}$.
    Plugging back in to \cref{eq:sup_of_eps} gives
    \[
        \chi(m;\vec{t}) < 3 +  c \pi \sup_{\tau > 0} \tau \eps(\tau;\vec{t})  < 3+4\sqrt{d}
        \; ,
    \]
    as claimed.
    
    It remains to prove that $\chi(m;\vec{t}) < 3 + 8\sqrt{m} \log'(d/m)$. Notice that this bound is weaker than the above if $d \leq 4 m\log'(d/m)^2$. So, in what follows, we may assume that $d >4  m\log'(d/m)^2$. This in particular implies that, say, $d > 70m$, which implies that $\tau^* = \log(2(d-m)/m)/\pi > 1.5 > 1$ and $q^* := e^{-\pi \tau^*} = m/(2(d-m)) < m/(1.9d)$. We may also assume without loss of generality that $\vec{t} \in [0,1/2]^d$ so that $\min_{\vec{z} \in \Z^d} \|\vec{z} - \vec{t}\|^2 = \|\vec{t}\|^2$. Since we have already assumed that $\min_{\vec{z} \in \Z^d} \|\vec{z} - \vec{t}\|^2 < m$, it follows that $m > \|\vec{t}\|^2$.

    With all of this out of the way, we see that in this case, we may restrict our attention to $\tau > \tau^* > 1$ so that $q := e^{-\pi \tau} < q^* < m/(1.9d)$. We may then apply \cref{lem:var_bound_big_tau_primal}, which yields
    \begin{align*}
       c\pi \tau \eps(\tau;\vec{t}) &<  2\pi c \tau \sqrt{4\|\vec{t}\|^2/\tau^2 + 6dq} \\
       &< 2\pi c \sqrt{4m + 6 \tau^2 dq} \\
       &< 2\pi c \sqrt{4m + 6 (\tau^*)^2 dq^*} \\
       &< 2\pi c \sqrt{4m + 6 (\tau^*)^2 m/1.9} \\
       &< 2\pi c \tau^* \cdot \sqrt{4m/1.5^2 + 6m/1.9} \\
        &< 4.6 \log(2d/m) \sqrt{m}\\
        &<8\sqrt{m} \log(d/m)
        \; .
    \end{align*}
    Plugging back in to \cref{eq:sup_of_eps} yields 
    \[
        \chi(m;\vec{t}) < 3 + 8 \sqrt{m} \log(d/m)
        \; ,
    \]
    as needed.
\end{proof}

\section{Some plots}
\label{sec:some_plots}

\begin{figure}
    \centering
    \includegraphics[width=\textwidth]{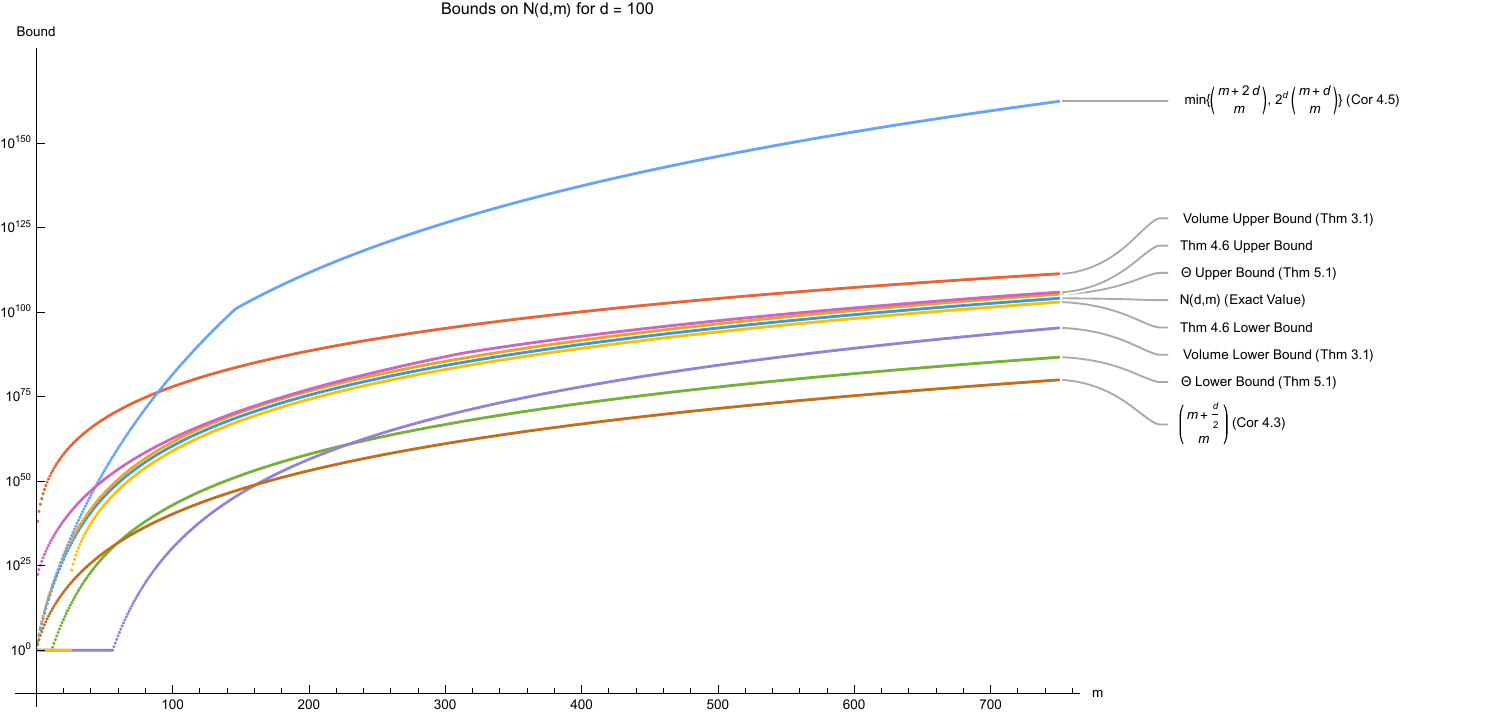}
    \caption{A plot of all bounds on $N(d,m)$ in this paper for dimension $d = 100$. See \cref{fig:d100-largem} for a similar plot for much larger $m$.
    (In all of our plots, we assume that lower bounds are always at least one. When our bounds are on $r_d(m)$ rather than $N(d,m)$, we obtain the bound on $N(d,m) $ by summing over all $m' \leq m$.)}
    \label{fig:d100-smallm}
\end{figure}

\begin{figure}
    \centering
    \includegraphics[width=\textwidth]{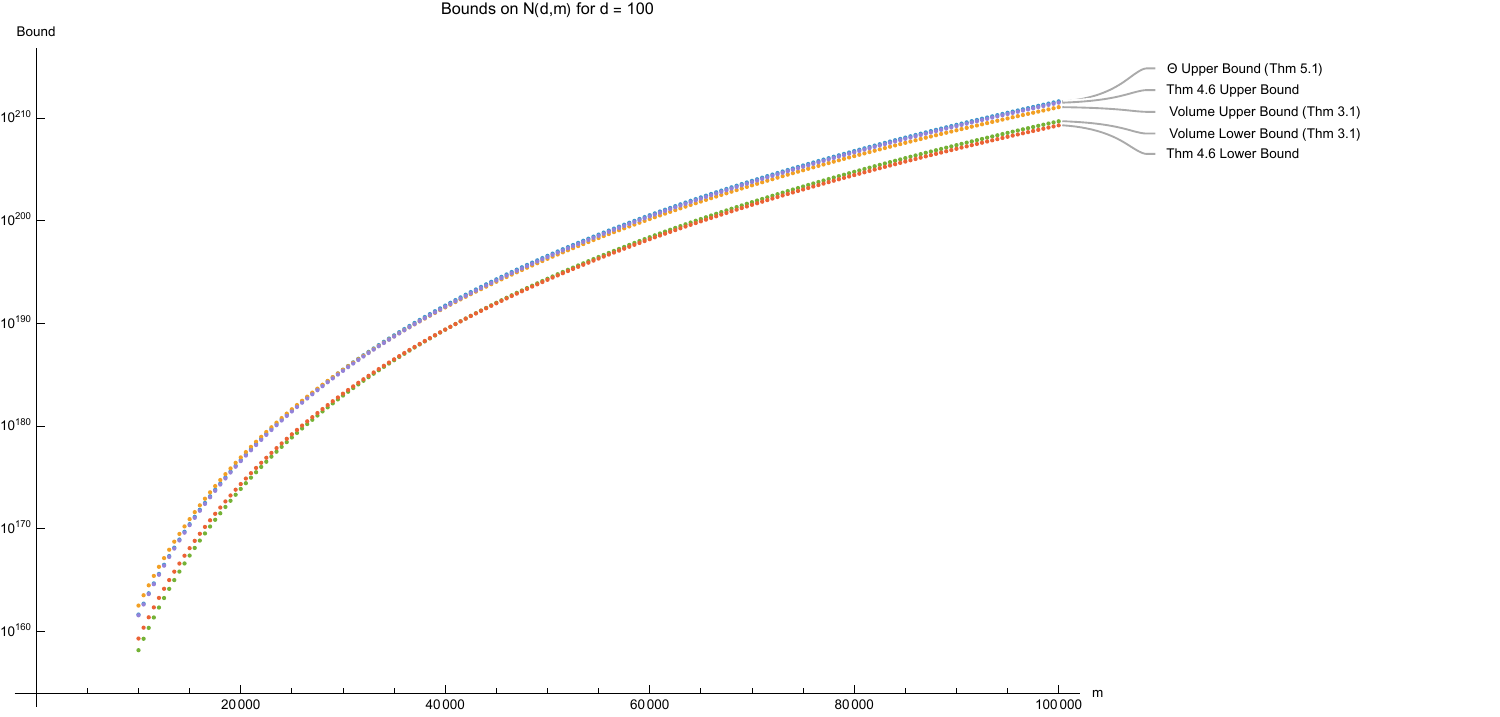} 
    \caption{A plot of bounds on $N(d,m)$ for $d = 100$ showing the behavior as $m$ becomes quite large relative to $d$. Eventually the volume-based bounds become quite tight. See \cref{fig:d100-smallm} for a similar plot for much smaller values of $m$.
    (We have elided all but the tightest bounds for these particular parameters, as well as the exact value of $N(d,m)$, which is already rather difficult to compute at this scale.)}
    \label{fig:d100-largem}
\end{figure}

\begin{figure}
    \centering
    \includegraphics[width=\textwidth]{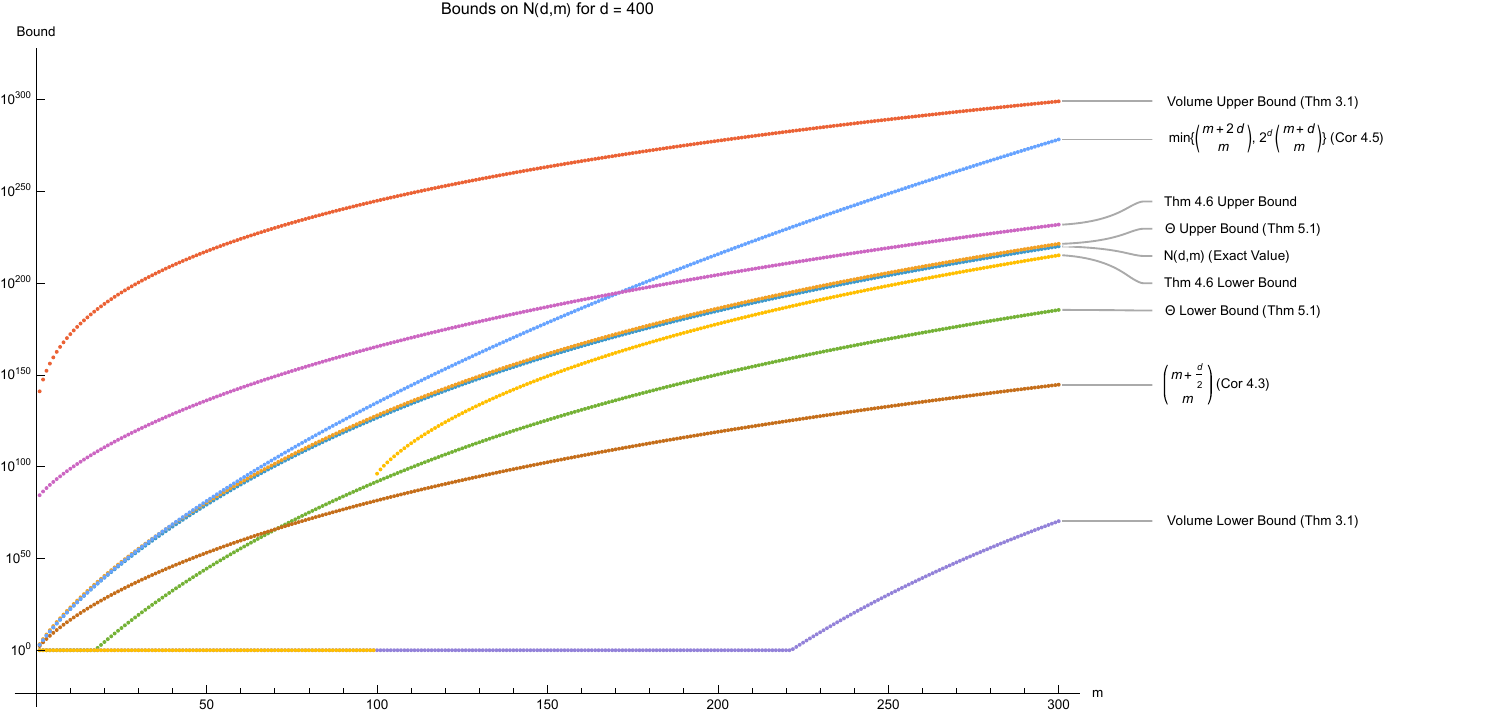} \\
    \includegraphics[width=\textwidth]{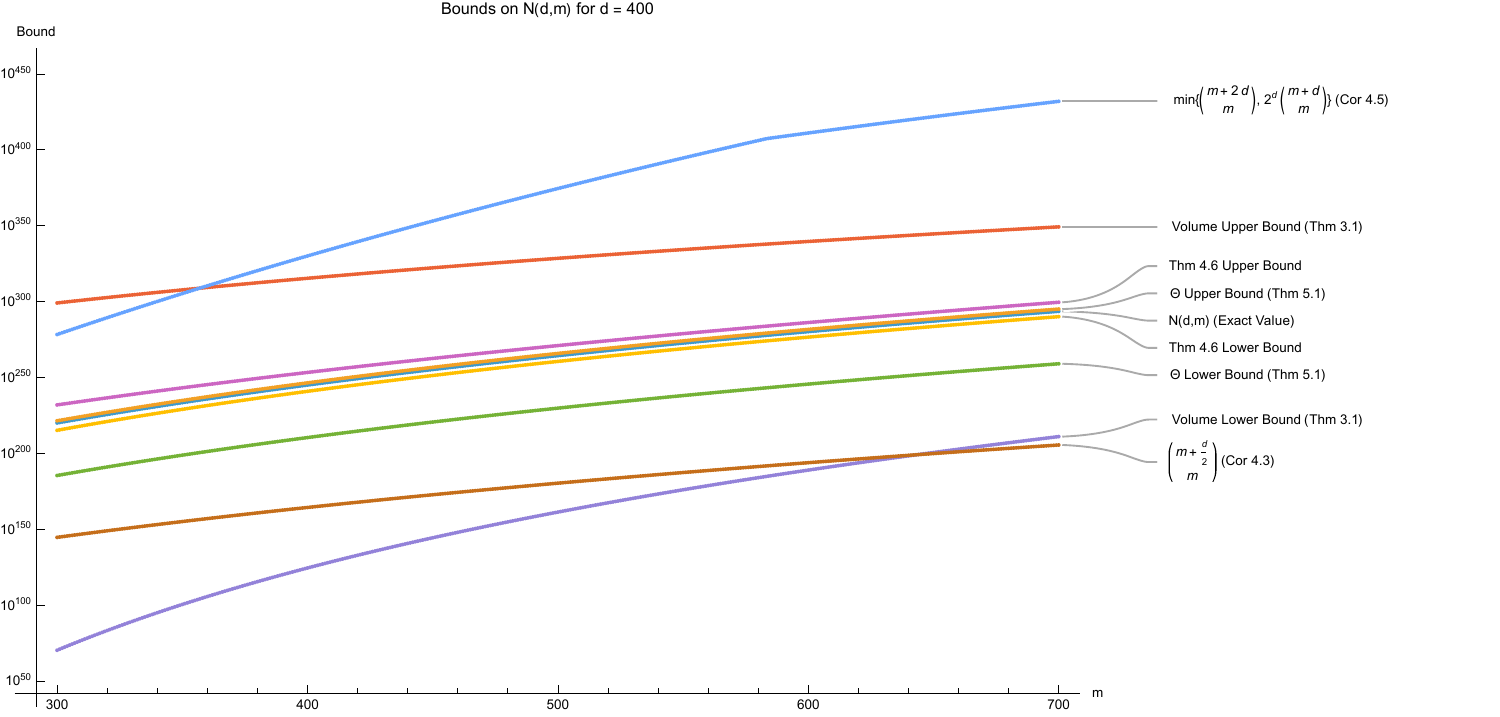} 
    \caption{Two plots showing the behavior of all of the bounds on $N(d,m)$ in this paper for $d = 400$.}
    \label{fig:d400}
\end{figure}

\begin{figure}
    \centering
    \includegraphics[width=\textwidth]{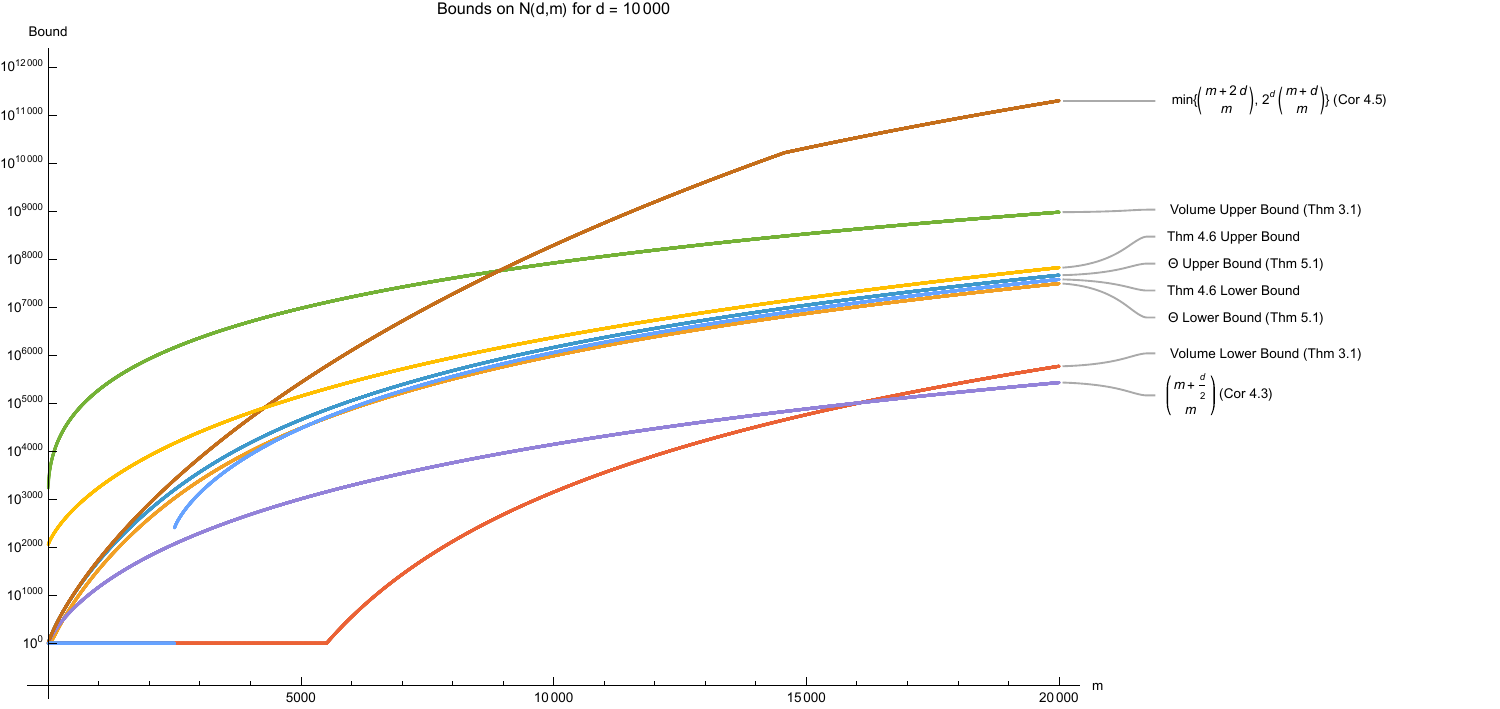}
    \caption{A plot showing the behavior of all of the bounds on $N(d,m)$ in this paper for $d = 10000$.}
    \label{fig:d10000}
\end{figure}

\begin{figure}
    \centering
    \includegraphics[width=\textwidth]{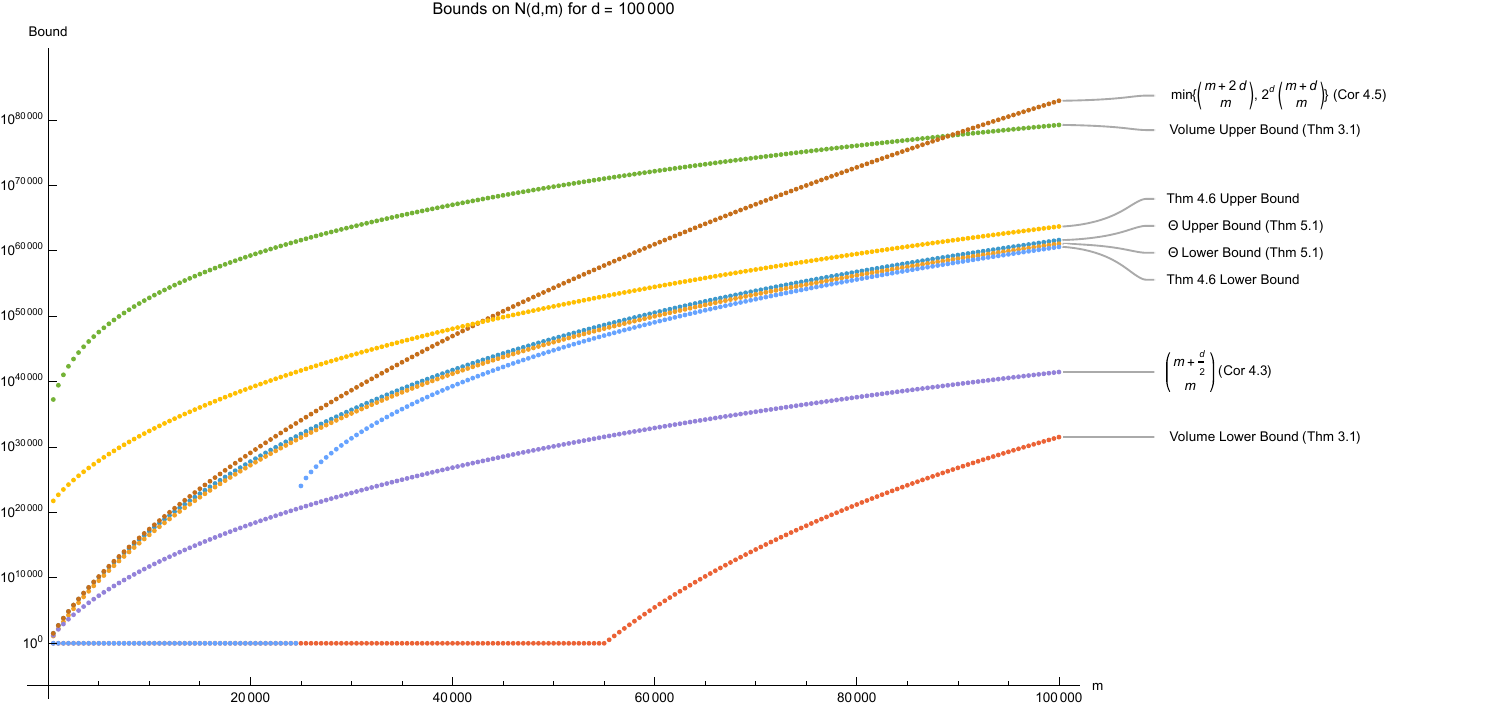}
    \caption{A plot showing the behavior of all of the bounds on $N(d,m)$ in this paper for $d = 100000$. See also \cref{fig:d100000-largem} for a similar plot for larger values of $m$.}
    \label{fig:d100000-smallm}
\end{figure}

\begin{figure}
    \centering
    \includegraphics[width=\textwidth]{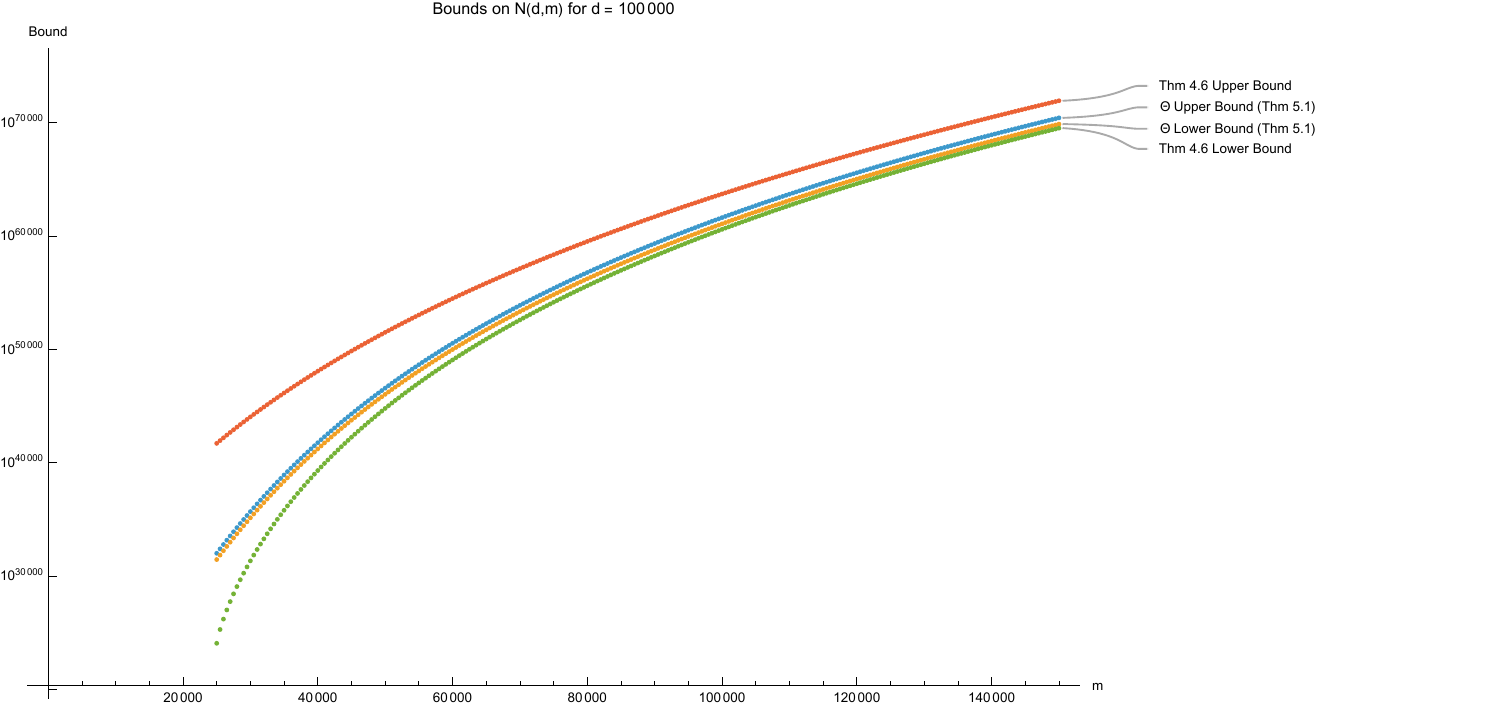}
    \caption{A plot showing the behavior of bounds on $N(d,m)$ for $d = 100000$. See also \cref{fig:d100000-smallm} for a similar plot for smaller values of $m$.}
    \label{fig:d100000-largem}
\end{figure}

Finally, we collect some plots showing how our bounds on $N(d,m)$ compare in different settings. These plots, shown in \cref{fig:d100-smallm,fig:d100-largem,fig:d400,fig:d10000,fig:d100000-smallm,fig:d100000-largem}, provide a nice intuitive understanding of when different bounds are useful. (At the risk of overwhelming the reader, we have endeavored to include a variety of plots to reflect the variety of bounds that we present and the variety of use cases that we anticipate.)

The plots (of course) show what one would expect. In particular, the very simple and well-known volume-based bounds from \cref{thm:volume_estimate} are quite loose until $m$ becomes much larger than $d$. (See \cref{fig:d100-largem} for a nice illustration of when these bounds start to become stronger than our other bounds.)  In contrast, our simple bounds in terms of binomial coefficients from \cref{cor:simple_binomial_lower,cor:simple_binomial_upper} are reasonably tight when $m$ is small but rather loose when $m$ is large. (The main benefit of these bounds is, of course, their simplicity.) 

The more complicated binomial-based bounds from \cref{thm:exponential_binomial} are quite tight. Asymptotically, it is easy to see that they are loose by a factor of $2^{\Theta(d)}$ . But, the hidden constant in this error term is small enough that these bounds are \emph{quite} tight in practice even for quite large parameters, and they are in fact our tightest bounds for a wide range of parameters. Of course, this is due to the fact that these bounds follow from carefully applying exact formulas for small values of $d$. (The lower bound is, however, zero for $m \lesssim d/4$.)

The theta-function-based bounds following Mazo and Odlyzko from \cref{thm:main_theta_theorem} are very tight asymptotically, with error less than $2^{O(\min\{\sqrt{d}, \sqrt{m} \log d \})}$. Indeed, the upper bound from  \cref{thm:main_theta_theorem} (which is more-or-less trivial to prove) is remarkably tight in essentially all parameter regimes. The lower bound is not as tight as the binomial-based lower bound from \cref{thm:exponential_binomial} until $d$ becomes quite large. In particular, \cref{fig:d100000-smallm,fig:d100000-largem} show that for $d = 100000$, the upper and lower bounds are really quite tight (as one can of course deduce directly from \cref{thm:main_theta_theorem}).

We should emphasize that, though we are in fact comparing these bounds, they are in many ways incomparable. In particular, the bounds from \cref{thm:volume_estimate,thm:main_theta_theorem} apply to the number of lattice points $N(d,m;\vec{t})$ in a ball centered around an arbitrary point $\vec{t}$, while the binomial-based bounds from \cref{sec:binomials} only apply to the centered case when $\vec{t} = \vec0$, but also can be used to estimate the number $r_d(m)$ of lattice points on a \emph{sphere}. Furthermore, the bounds evidently vary in their simplicity and ease of use.

\appendix

\crefalias{section}{appendix}
\crefalias{subsection}{appendix}

\section{Additional proofs}

    \subsection{Proof of \texorpdfstring{\cref{lem:r_7_and_r_9_lower_bound} (lower bounds on $r_7$ and $r_9$)}{lower bounds in 7 and 9 dimensions}}
    \label{app:r_7_and_r_9_lower_bounds}
    
    For both the $d=7$ case and the $d=9$ case, we will follow the methods of \cref{r_5_bound}. 
    
    First we will set $d=7$ and we will show the bound holds for $m \geq 100$ as direct computation (e.g., using \cref{eq:r_d_recurrence}) shows the bound holds for $m < 100$.
    
    Let $k := \floor{\sqrt{m}}$. Then, we can apply \cref{eq:r_d_recurrence} and use the lower bound for $r_6(m)$ from \cref{binomial_r_6_r_8_bounds} to get 
    \begin{align*}
        r_7(m) &= \sum_{i=0}^{m}{r_1(i) \cdot r_6(m-i)} \\
        &=r_6(m) + 2 \sum_{j=1}^{k}{r_6(m-j^2)} \\
        &\geq {m+2 \choose 2} + 2\sum_{j=1}^{k}{{m-j^2+2 \choose 2}} \\
        &= \frac{1}{15} \cdot (3 k^5+ 15 m^2k -10 mk^3) + \frac{(m-k^2)^2}{2} + \frac{8mk - 2k^3}{3} + \frac{3m-3k^2}{2} + \frac{44k}{30} + 1 \\
        &\geq \frac{1}{15} \cdot (3 k^5+ 15 m^2k -10 mk^3)\\
        &\geq \frac{1}{15} \cdot ( 8 m^{5/2} - 20m + 15\sqrt{m}-3)\\
        &\geq \frac{1}{2} m^{5/2}
        \;, 
    \end{align*}
    where the second-to-last inequality uses the fact that $k \geq \sqrt{m}-1$ and the last inequality holds for $m \geq 100$. 

    Finally, using the simple inequality $\binom{x}{3/2} \leq x^{3/2}$ for $x \geq 0$, we have
     \begin{align*}
        {m+5/2 \choose 5/2} &= \frac{2}{5} (m+5/2) \cdot  {m+3/2 \choose 3/2} \\
        &\leq \frac25 (m+5/2)(m+3/2)^{3/2} \\
        &\leq \frac{1}{2}m^{5/2}\; ,
    \end{align*}
    where the last inequality holds for $m \geq 100$. 
    The result follows.

    Analogously, we apply \cref{eq:r_d_recurrence} and the lower bound for $r_8(m)$ in \cref{r_8 bound} which yields
    \begin{align*}
        r_9(m) &= \sum_{i=0}^{m}{r_1(i) \cdot r_8(m-i)} \\
        &=r_8(m) + 2 \sum_{j=1}^{k}{r_8(m-j^2)} \\
        &\geq 14m^3 + 2 \sum_{j=1}^{k}{14(m-j^2)^3} \\
        &= 14m^3+\frac{2}{15}\left(-5k+35k^3-105k^5-105k^6-30k^7-21km+210k^3m+315k^4m+126k^5m\right) \\
        &\quad + \frac{2}{15}\left(-105km^2-315k^2m^2-210k^3m^2+210km^3 \right) \\
        &\geq 14m^3 + \frac{2}{15}\left(126 (\sqrt{m}-1)^5m - 30m^{7/2} - 420m^3 - 210m^{5/2}-21 m^{3/2}-5\sqrt{m} \right) \\ 
        &> \frac12 m^{7/2} 
        \; ,
    \end{align*}
    where we apply the fact that $\sqrt{m} - 1 \leq k \leq \sqrt{m}$ (and throw away some lower-order positive terms together with the clearly non-negative term $210km^3 - 210k^3 m^2$) for the penultimate inequality and the last inequality holds for $m \geq 100$.

    Once again, we can apply the simple inequality $\binom{x}{5/2} \leq x^{5/2}$ for, say, $x \geq 5/2$, we have
     \begin{align*}
        {m+7/2 \choose 7/2} &= \frac{2}{7} (m+7/2) \cdot  {m+5/2 \choose 5/2} \\
        &\leq \frac27 (m+7/2)(m+5/2)^{5/2} \\
        &\leq \frac{1}{2}m^{7/2}\; ,
    \end{align*}
    where the last inequality holds for $m \geq 100$. 
    Thus, the $d=9$ case is resolved as well.

    \subsection{Proof \texorpdfstring{of \cref{lem:finite_sums_primal_new,cor:k=1}}{that certain series can be approximated by partial sums}}
    \label{app:finite_sums_primal}

    \begin{proof}[Proof of \cref{lem:finite_sums_primal_new}]
    We have
    \[
        \Theta(\tau;t) - \Theta_k(\tau;t) = \sum_{z=k+1}^\infty q^{(z-t)^2} + \sum_{z=k}^\infty q^{(z+t)^2} = q^{(k+t)^2} + \sum_{z=k+1}^\infty \big( q^{(z-t)^2} + q^{(z+t)^2}\big)
        \; ,
    \]
    so that evidently this difference is non-negative. Bounding the above, we see that
    \[
         \Theta(\tau;t) - \Theta_k(\tau;t) \leq q^{k^2} + 2 \sum_{z=k+1}^\infty q^{(z-t)^2} \leq q^{k^2} + \frac{q^{w^2}}{\pi \tau w}
         \; ,
    \]
    as claimed,
    where the last inequality is \cref{clm:Gaussian_sum_bounds}.

    Now, let $g_\tau(x) := x^2 q^{x^2}$.  Similarly to the above, we have
    \[
        \mu(\tau;t) -  \mu_k(\tau;t) = g_\tau(k+t) + \sum_{z=k+1}^\infty \big( g_\tau(z-t) + g_\tau(z+t) \big)
        \; ,
    \]
    which is evidently non-negative. And, notice that $g_\tau(x)$ is a decreasing function of $x$ for $x \geq 1/\sqrt{\pi \tau}$. It follows that $g_\tau(k+t) \leq g_\tau(k)$ and for $z \geq 1 +t+ 1/\sqrt{\pi \tau}$, $g_\tau(z+t) + g_\tau(z-t)  \leq 2 g_\tau(z-t)$. Since $k \geq t + 1/\sqrt{\pi \tau}$, it follows that
    \[
        \mu(\tau;t) -  \mu_k(\tau;t) \leq g_\tau(k) + 2\sum_{z=k+1}^\infty g_\tau(z-t)
        \leq g_\tau(k) + \Big(1  + \frac{1}{2 \pi w^2 \tau}\Big) \cdot \frac{w q^{w^2}}{\pi \tau} \leq g_\tau(k) + \frac{5w}{4\pi \tau} \cdot q^{w^2}
        \; ,
    \]
    as needed,
    where the second inequality is \cref{clm:Gaussian_sum_bounds} and the last inequality uses the assumption that $w \geq \sqrt{2/(\pi \tau)}$. 

    Finally, let $f_\tau(x) := x^4 q^{x^2}$, and similarly notice that 
    \[
        M(\tau;t) -  M_k(\tau;t) = f_\tau(k+t) + \sum_{z=k+1}^\infty \big( f_\tau(z-t) + f_\tau(z+t) \big)
    \]
    is non-negative. $f_\tau(x)$ is a decreasing function of $x$ for $x \geq \sqrt{2/(\pi \tau)}$. So, $f_\tau(k+t) \leq f_\tau(k)$ and for $z \geq 1+t+\sqrt{2/(\pi \tau)}$, we have $f_\tau(z+t) + f_\tau(z-t) \leq 2f_\tau(z-t)$. Since $k \geq t + \sqrt{2/(\pi \tau)}$, it follows that
    \begin{align*}
        M(\tau;t) - M_k(\tau;t) &\leq f_\tau(k) + 2\sum_{z=k+1}^\infty f_\tau(z-t)\\
            &\leq f_\tau(k) + \Big( \frac{w^3}{\pi \tau} + \frac{3w}{2\pi^2 \tau^2} + \frac{3}{4 \pi^3 w \tau^3} \Big)  q^{w^2} \\
            &= f_\tau(k) +  \frac{w^3}{\pi \tau} \cdot \Big(1 + \frac{3}{2\pi w^2 \tau} + \frac{3}{4\pi^2 w^4 \tau^2} \Big)  q^{w^2}\\
            &\leq f_\tau(k) + \frac{31 w^3}{16\pi \tau}  q^{w^2}
        \; ,
    \end{align*}
    as claimed,
    where the second inequality is \cref{clm:Gaussian_sum_bounds} and the last inequality uses the assumption that $w \geq \sqrt{2/(\pi \tau)}$. 
    \end{proof}

\begin{proof}[Proof of \cref{cor:k=1}]
    (It would be slightly more natural to obtain bounds in terms of $q^{(1+t)^2}$ rather than $q^{1+t^2}$. However, we are careful to obtain bounds in terms of $q^{1+t^2}$ because (1) this is precisely what we need; and (2) this achieves much better constants, which yields tighter results when $t$ is small. We recognize that this unfortunately looks like a typo, but it is not.)

     Notice that
     \[
         0 \leq \Theta_2(\tau;t) - \Theta_1(\tau;t) = q^{(1+t)^2} +q^{(2-t)^2} 
        < 1.01 q^{1+t^2}
         \; ,
     \]
     where in the last step we have used the fact that the function 
     \[
         t \mapsto \frac{q^{(1+t)^2} + q^{(2-t)^2}}{q^{1+t^2}}
     \]
     is convex on the interval $t \in [0,1/2]$ so that it is maximized at either $t = 0$ or $t = 1/2$. One then checks directly that the bound holds for $t \in \{0,1/2\}$.
     Applying \cref{lem:finite_sums_primal_new} with $k = 2$ immediately implies that
     \[
         0 \leq \Theta(\tau;t) - \Theta_1(\tau;t) \leq q^4 + 1.01 q^{1+t^2} + \frac{1}{\pi \tau (2-t)} \cdot q^{(2-t)^2} < 1.1 q^{1+t^2}
         \; ,
     \]
     as claimed.

     Similarly, we have 
     \[
         0 \leq \mu_2(\tau;t) - \mu_1(\tau;t) =  (1+t)^2 q^{(1+t)^2} + (2-t)^2 q^{(2-t)^2} < 1.01 q^{1+t^2}
         \; ,
     \]
     where again the last step follows by convexity.
     Applying \cref{lem:finite_sums_primal_new} again yields
     \[
     0 
         \leq \mu(\tau;t) - \mu_1(\tau;t) 
         < 4 q^4 + 1.01 q^{1+t^2} + \frac{5(2-t)}{4 \pi \tau} q^{(2-t)^2}
         < 1.1 q^{1+t^2}
     \; ,
     \]
     as claimed.

     Finally, we have 
     \[
     0 \leq M_2(\tau;t) - M_1(\tau;t) = (1+t)^4 q^{(1+t)^2} +  (2-t)^4  q^{(2-t)^2} < 1.002 q^{1+t^2}
     \; ,
     \]
     where the last step again works by convexity (and where we have been slightly more careful with the constants here).
     \cref{lem:finite_sums_primal_new} then implies that
      \[
         0 \leq M(\tau;t) - M_1(\tau;t) < 1.002q^{1+t^2}   + 16 q^{4} + \frac{31 (2-t)^3}{16 \pi \tau} q^{(2-t)^2} < 1.1 q^{1+t^2}
         \; ,
      \]
      as needed.
\end{proof}

\subsection{Good approximations for small \texorpdfstring{$\tau$}{tau}}
\label{app:finite_sums_dual}

We now show an analogue of the finite approximations given in \cref{sec:finite_sum_defs} and \cref{lem:finite_sums_primal_new} that works in the case when $\tau$ is small. These work by applying the Poisson summation formula, as in \cref{fact:PSF_theta_mu_M}.

We define
\[
    \widehat{\Theta}_k(\tau;t) := \frac{1}{\sqrt{\tau}} + \frac{2}{\sqrt{\tau}} \cdot \sum_{z=1}^k \cos(2\pi t z) p^{z^2} \; , \quad \widehat{\mu}_k(\tau;t) := \frac{1}{2\pi \tau} \cdot \widehat{\Theta}_{k}(\tau;t) - \frac{2}{\tau^{5/2}} \cdot \sum_{z=1}^k \cos(2\pi z t)z^2 p^{z^2}
    \; ,
\]
and
\[
\widehat{M}_k(\tau;t) := \frac{3}{\pi \tau} \cdot \widehat{\mu}_k(\tau;t) - \frac{3}{4\pi^2 \tau^2} \cdot \widehat{\Theta}_{k}(\tau;t)  + \frac{2}{\tau^{9/2}} \sum_{z=1}^k \cos(2\pi z t) z^4 p^{z^2}
\; ,
\]
where $p := e^{-\pi/\tau}$. These expressions arise from taking Fourier transforms of the respective summands defining $\Theta$, $\mu$, and $M$, as in \cref{fact:PSF_theta_mu_M}. 

The following lemma shows that these finite sums yield good approximations to the relevant series when $k$ is not too small relative to $\sqrt{\tau}$ (so, in particular, these bounds are more useful when $\tau$ is smaller).

\begin{lemma}
    \label{lem:finite_sums_dual}
    For any $\tau > 0$, $t \in \R$, and integer $k \geq \sqrt{2\tau/\pi}$, we have
    \[
        \big|\Theta(\tau;t) - \widehat{\Theta}_k(\tau;t)\big| \leq \frac{\sqrt{\tau}  }{\pi k} \cdot p^{k^2}
        \; ,
    \]
     \[
    \big|\mu(\tau;t) -  \widehat{\mu}_k(\tau;t)\big| \leq \frac{3k }{2\pi \tau^{3/2}} \cdot p^{k^2}
        \; ,
    \]
    and
    \[
         \big|M(\tau;t) -  \widehat{M}_k(\tau;t)\big| \leq \frac{35 k^3 }{8 \pi \tau^{7/2}} \cdot p^{k^2} 
        \; ,
    \]
    where $p := e^{-\pi/\tau}$.
\end{lemma}
\begin{proof}
    By \cref{fact:PSF_theta_mu_M}, we can write
    \[
        \Theta(\tau;t) = \frac{1}{\sqrt{\tau}} + \frac{2}{\sqrt{\tau}} \sum_{z = 1}^\infty \cos(2\pi zt) p^{z^2}
        \; .
    \]
    Therefore,
    \[
        \big|\Theta(\tau;t) - \widehat{\Theta}_k(\tau;t)\big| \leq \frac{2}{\sqrt{\tau}} \cdot \Big|\sum_{z=k+1}^\infty \cos(2\pi zt) p^{z^2}\Big| \leq \frac{2}{\sqrt{\tau}} \cdot  \sum_{z=k+1}^\infty p^{z^2} \leq \frac{\sqrt{\tau}}{\pi k} \cdot p^{k^2}
        \; ,
    \]
    as needed, where the last inequality is \cref{clm:Gaussian_sum_bounds}.

    Now, let $g_\tau(x) := x^2 q^{x^2}$. By \cref{fact:PSF_theta_mu_M}, we have
    \[
        \mu(\tau;t) = \frac{\Theta(\tau;t)}{2\pi \tau} - \frac{2}{\tau^{5/2}} \cdot \sum_{z=1}^\infty \cos(2\pi zt) g_{1/\tau}(z)
        \; .
    \]
    Therefore,
    \begin{align*}
        |\mu(\tau;t) - \widehat{\mu}_k(\tau;t)| 
            &\leq \frac{1}{2\pi \tau} \cdot \big|\Theta(\tau;t) - \widehat{\Theta}_k(\tau;t)\big| +  \frac{2}{\tau^{5/2}} \cdot \Big| \sum_{z=k+1}^\infty \cos(2\pi z t) g_{1/\tau}(z) \Big| \\
                &\leq \frac{p^{k^2}}{2\pi^2 k \tau^{1/2}} + \frac{2}{\tau^{5/2}} \cdot \sum_{z=k+1}^\infty g_{1/\tau}(z)\\
                &\leq \frac{p^{k^2}}{2\pi^2 k \tau^{1/2}} + \frac{k p^{k^2}}{\pi \tau^{3/2}} + \frac{p^{k^2}}{2\pi^2 k\tau^{1/2}}\\
                &\leq 
                \frac{3k p^{k^2}}{2\pi \tau^{3/2}}
        \; ,
    \end{align*}
    where the second inequality is \cref{clm:Gaussian_sum_bounds} and the last inequality uses the fact that $k \geq \sqrt{2\tau/\pi}$.

    Finally, let $f_\tau(x) := x^4 q^{x^2}$. By \cref{fact:PSF_theta_mu_M}, we have
    \[
        M(\tau;t) = \frac{3}{\pi \tau} \cdot \mu(\tau;t) - \frac{3}{4\pi^2 \tau^2} \cdot \Theta(\tau;t)+ \frac{2}{\tau^{9/2}} \cdot \sum_{z=1}^\infty \cos(2\pi z t) f_{1/\tau}(z)
        \; .
    \]
    Therefore,
    \begin{align*}
        |M(\tau;t) - \widehat{M}_k(\tau;t)| 
            &\leq \frac{3}{\pi \tau} \cdot \big|\mu(\tau;t) - \widehat{\mu}_k(\tau;t)\big| + \frac{3}{4\pi^2 \tau^2} \cdot \big|\Theta(\tau;t) - \widehat{\Theta}_k(\tau;t)\big| + \frac{2}{\tau^{9/2}} \cdot \sum_{z=k+1}^\infty  f_{1/\tau}(z)\\
            &\leq  \frac{9k p^{k^2}}{2\pi^2 \tau^{5/2}} +  \frac{3 p^{k^2} }{4\pi^3 k \tau^{3/2}}   + \frac{2}{\tau^{9/2}} \cdot \sum_{z=k+1}^\infty f_{1/\tau}(z)\\
            &\leq \frac{9k p^{k^2}}{2\pi^2 \tau^{5/2}} +  \frac{3 p^{k^2} }{4\pi^3 k \tau^{3/2}} + \Big(  \frac{k^3}{\pi \tau^{7/2}} + \frac{3k}{2\pi^2 \tau^{5/2}} + \frac{3}{4\pi^3 k \tau^{3/2} }  \Big) \cdot p^{k^2}\\
            &\leq \frac{35 k^3 p^{k^2}}{8 \pi \tau^{7/2}}
        \; ,
    \end{align*}
    as needed,
    where the second inequality is \cref{clm:Gaussian_sum_bounds} and the last inequality uses the assumption that $k \geq \sqrt{2\tau/\pi}$.
\end{proof}

\end{document}